\documentclass[a4paper,11pt,leqno]{article}

\usepackage{amsmath} % omitted
\usepackage{amssymb} % omitted
\usepackage{amscd} % omitted
\usepackage{latexsym} % omitted

\usepackage{array} % omitted

\usepackage{mathrsfs}  % omitted

\usepackage[utf8]{inputenc}

\numberwithin{equation}{section}

\numberwithin{equation}{section}  

\allowdisplaybreaks

\title{Lineage of the Theory of Invariant Integrals on 
Groups\\  
{\large  
Hurwitz, Schur, Weyl, Haar, Neumann, Kakutani, Weil, Kakutani-Kodaira}
}
\author{Takeshi HIRAI}

\date*{}

\begin{document}

\maketitle

\setcounter{section}{0}
\setcounter{page}{1}

{\bf Abstract.} 
This is mainly a translation of Proc.\;of the 29th Symp.\;on History of Mathematics, Tsuda University, held Oct.\;2018, of my talk.\footnote{Cf. \cite{Hir18} in {\bf References}.}  
The first occasion when I studied the history of the theory of invariant integrals (or measures) was an unintended opportunity where I was asked to write \lq\lq\,Kaisetsu\,''\hspace{.2ex}(explanatory and commentary article) 
to the new book\footnote{
This book \cite{Sai15} was published as a volume of the series of pocket-size book of reproduction of natural science classics, called {\it Chikuma Gakugei Bunko}. Each book has an additional part called {\it Kaisetsu} of about 20 pages, written by an expert of the subject of the book and attracts readers in addition to those of the main body.  In this case I wrote an article of 
30 pages and Prof.\;Saito thanked me writing \lq I think it's a very interesting commentary with rich of contents' in his afterword. 
To prepare this article, I read several mathematical papers, in particular around Weil representations, and 
also Weil's {\it Souvenirs d'apprentissage},\, {\it Commentaire}\, in his {\it Collected Papers} Vols I\,--\,III 
and so on. 
Moreover, Mr.\;I.\;Ebihara, the editor of the book, wrote in the back-page-campaign of the bookcover as \lq{\it The study of integrals on the group-space has begun with Hurwitz 
at the end of 19th Century. $\cdot\,\cdot\,\cdot$ }', asking me as \lq Is that correct ?\,' I answered \lq\,Yes\,\rq, not knowing well about his work. And so naturally I was forced later to study his paper [Hur1897] to ensure my answer \lq Yes\rq. 
} 
of Prof.\;M.\;Saito, a first translation into Japanese of the famous Weil's book \lq{\it L'int\'egrations dans les groupes topologiques et $\cdots\cdots$}'. Personally, for my professional work, it was needed only to read this Weil's original book and several text books on measure theory. For writing the Kaisetsu mentioned above, other than several mathematical papers 
and also Weil's non-mathematical works, 
it was sufficient for me to 
read Haar's original paper roughly and similarly for other historical classics. Thus I felt the necessity of further study and now  
 I make intendedly a historical study. Reading original papers due to Hurwitz, Schur, Weyl and so on in detail, I explain their contents from my point of view, and give the relationships among them.  

\vskip2em
\newpage

{\small 

\tableofcontents

}

\section{Hurwitz's pioneering work on integrals on Lie groups}

The title of his paper [Hur1897] is \lq\,{\it On the generation of invariants through integration}\,' and its main objective is to generate invariant polynomials under natural actions of Lie groups $G$, especially rotation groups $SO(n)$ and also unitary groups $SU(n)$.

For this end, he proposed to use invariant integrals on Lie groups $G$.

\vskip1.2em
This 20 pages paper has an \lq\,Introduction'  of about 1 page and is not separated by sections. But it has several separations indicated by a horizontal segment in the middle of blank lines. Take appropriate 4 separations as sections and put appropriate titles to each sections, then the structure of the paper is approximately as follows (the technical terms are modernized as those of today):
\vskip1em 

$
\begin{array}{lll}
&{\rm Introduction} & \mbox{\rm p.70}
\\
\S\,1&\mbox{\rm Generation of invariants on continuous groups using invariant} 
&
\\
\hspace*{5ex} &\mbox{\rm integration, in particular, case of $SO(n)$}& 
\mbox{\rm p.71}
\\ 
\S\,2&\mbox{\rm Global coordinates on $SO(n)$ and formula of invariant measure}
\;& \mbox{\rm p.75}
\\
\S\,3&\mbox{\rm Global coordinates on $SU(n)$ and formula of invariant measure}
\;& \mbox{\rm p.80}
\\
\S\,4&\mbox{\rm General theory for Lie groups}& \mbox{\rm p.86}
\end{array}
$

\setcounter{subsection}{-1}
\subsection{On \lq\,Introduction'}

Hurwitz explained first as 
\begin{quotation} 
Assume that a finite group $G$ acts on variables $x_1,x_2,\ldots,x_n$ by permuting them. Then $G$ acts on a function $f(x_1,x_2,\ldots,x_n)$, and we get an invariant function $\varphi$ by averaging $gf$ over $g\in G$. This process can be applicable in some cases of infinite number of variables.
\end{quotation}  
Then continued, in German original, as (use machine translation if necessary)
\begin{quotation} 
\noindent
$\cdots\cdots$ Ich habe nun Gedanke verfolgt, dieses sich so zu sagen von selbst darbietende Verfahlen zur Erzeugung der Invarianten auf die continuierlichen Gruppe zu \"ubertragen, wo dann naturgem\"a\ss\/ bestimmte Integrale an die Stelle der Summen treten. 
Dabei richtete ich mein Augenmerk zun\"achst 
auf die ganzen rationalen Invarianten der algebraischen Formen, 
also auf diejenigen ganzen rationalen Functionen der Coefficienten einer Form, 
die sich nicht \"andern, wenn man auf die Variabeln der Form eine beliebige lineare, homogene, unimodulare Substitution aus\"ubt. 
Die Untersuchung f\"uhrte mich in dessen bald dazu,\, $\cdots\cdots$
\end{quotation} 

\noindent
then continued, in short, as 
\begin{quotation}
\noindent
More extended cases we study are in particular rotation groups $SO(n)$, 
 and also special unitary groups $SU(n)$. 
\end{quotation}

\subsection{On \lq\,\S 1 Generation of invariants on continuous groups using invariant integration$\cdots$'} 

Please admit me to use modern notation $G=SO(n)$. 
Take an element $g=(r_{\alpha,\beta})_{\alpha,\beta\in{\boldsymbol I}_n}$ of the rotation group $G$, with ${\boldsymbol I}_n:=\{1,2,\ldots,n\},$ and consider its natural action on 
the variables $x:={}^t(x_1,x_2,\ldots,x_n)$ under  
\begin{gather}
\label{2021-06-21-1}
x_\alpha=\sum_\beta r_{\alpha\beta}x'_\beta\quad\;{\rm or}\;\quad x=gx', 
\\
\label{2021-06-21-2}
 \sum_\gamma r_{\alpha\gamma} r_{\beta\gamma}=\delta_{\alpha\beta}\,,
\quad  |r_{\alpha\beta}|=1, 
\end{gather}
where $|r_{\alpha\beta}|$ denotes determinant. 
The object (\lq\lq\,Gebilde\,\rq\rq) in ${\boldsymbol R}^{n^2}$ of $(r_{\alpha\beta})_{\alpha,\beta\in{\boldsymbol I}_n}$ defined by the fundamental equation (\ref{2021-06-21-2}) is denoted by $R$, which is a submanifold corresponding to $G$. 
A measure on $R$ is denoted by $dR$ in \cite{Hur1897}, \,but permit me to denote it also as $dg\;(g\in G)$.  When the substitution (\lq\lq\,Substitution\,\rq\rq)\; 
(\ref{2021-06-21-1}) is put together with  
\begin{gather}
\label{2021-08-10-2}
x'_\alpha=\sum_\beta s_{\alpha\beta}x''_\beta\,,\qquad h=(s_{\alpha\beta})
\end{gather}
we have $gh=(r'_{\alpha\beta})$ with 
\begin{gather}
\label{2021-08-10-3}
x_\alpha=\sum_\beta r'_{\alpha\beta}x''_\beta\,,\qquad 
r'_{\alpha\beta}=\sum_\gamma r_{\alpha\gamma}s_{\gamma\beta}\,.
\end{gather}

Now consider a form 
\begin{eqnarray}
\label{2021-08-10-1}
\Phi(a;x)=\sum_{1\leqslant k\leqslant m}a_j\,p_j(x)
\end{eqnarray}
of polynomials in $x={}^t(x_1,x_2,\ldots,x_n)$ of homogeneous degree $p\ge 0$ with coefficients $a={}^t(a_1,a_2,\ldots,a_m)$, where  
$$
m=\dfrac{(n+p-1)!}{p!(n-1)!},\qquad 
p_j(x)\;\;\mbox{\rm different monomials of degree $p$}, 
$$

When $x$ is replaced by $x'$ under (\ref{2021-06-21-1}), we have 
\begin{eqnarray}
\nonumber
p_j(x)
&\to& 
p\Big({\sum}_\gamma r_{1\beta}x'_\beta,\,{\sum}_\gamma r_{2\beta}x'_\beta,\,\ldots,\,{\sum}_\gamma r_{n\beta}x'_\beta\Big)
\\
\nonumber
&& =\sum_k P_g^{jk}\,p_k(x'),\qquad P_g:=(P_g^{jk})_{j,k\in{\boldsymbol I}_m},
\\
\label{2021-08-10-4} 
\therefore
&& 
\Phi(a;x)=\Phi(a';x'),\qquad a'={}^tP_g^{\,-1}a.
\end{eqnarray}
Here we have\; $P_gP_h=P_{gh}$\,, and every element of $P_g$ is a polynomial in $r_{\alpha\beta}$\,.

Now let $F(a)$ be an arbitrary polynomial in $a={}^t(a_1,a_2,\ldots,a_m)$ of coefficients in $\Phi(a;x)$. Then, under \lq\lq\,Substitution\,\rq\rq\,:\; $x\to x'$ in (\ref{2021-06-21-1}), $F(a)$ is transformed to $F(a')$. The first assertion of Hurwitz is written in p.74, with an invariant measure $dR$ on $R$, as 
\begin{quotation}
D\,a\,s \"u\,b\,e\,r d\,a\,s G\,e\,b\,i\,l\,d\,e $R$ a\,u\,s\,g\,e\,d\,e\,h\,t\,e I\,n\,t\,e\,g\,r\,a\,l 
\begin{eqnarray}
\label{2021-08-11-1}
J(a)=\int F(a')\,dR 
\end{eqnarray}
s\,t\,e\,l\,l\,t n\,u\,n e\,i\,n\,e o\,r\,t\,o\,g\,o\,n\,a\,l\,e I\,n\,v\,a\,r\,\,i\,a\,n\,t\,e d\,e\,r F\,o\,r\,m $\Phi(a;x)$ d\,a\,r.
\end{quotation} 
Every \lq\lq\,orthogonale Invariante der 
Form $\Phi(a;x)$\,\rq\rq\/ is given by this integration. 

But \lq\lq das Abz\"ahlungsproblem der Anzahl der orthogonalen Invarianten der Form $\Phi(a;x)$\rq\rq\/ cannot be solved at this stage, whereas the finiteness of this number is secured by the first fundamental theorem of Hilbert. 
\vskip1em

{\bf Explanation 1.1.}\; Let $\pi$ be a natural matrix representation $\pi(g):=g\;(g\in G)$ on the representation space $V(\pi):={\boldsymbol C}^n$. Then the space of polynomials of degree $r$ of variables $x_1,x_2,\ldots,x_n$ is nothing but the space of symmetric tensor product $S^rV(\pi)$ in $\otimes^rV(\pi)$. The form $\Phi(a;x)$ gives the dual pairing of $S^rV(\pi)$ and its dual $\big(S^rV(\pi)\big)^*$, and both spaces receive the natural action of $G=SO(n)$, and the actions of both sides are denoted by $T(g)$ and $T^\vee(g):={}^tT(g)^{-1}$ respectively.  The replacement $F(a)$ by $F(a')$ can be expressed as $(g^{-1}F)(a):=F(a')=F\big(T^\vee(g)a\big)$, and the integration (\ref{2021-08-11-1}) is expressed as 
\begin{gather}
\label{2021-08-11-2}
J(a)=\int_G (g^{-1}F)(a)\,dg\quad\;{\rm or}\;\quad 
J=\int_G g^{-1}F\,dg \;\;\Big(\!=\int_G gF\,dg\Big). 
\end{gather} 

Das Abz\"ahlungsproblem f\"ur \lq\lq\,orthogonale Invariante der 
Form $\Phi(a;x)$\,\rq\rq\/ is proposed.

\subsection{On \lq\,\S 2\; Global coordinates on $SO(n)$ and formula of invariant measure\,\rq}

In general, let $\xi=(\xi_1,\xi_2,\ldots,\xi_\tau)$ be orthonormal coordinates in $\tau$-dimensional ${\boldsymbol R}^\tau$ and consider a submanifold ${\mathfrak{R}}$ of dimension $\sigma$ defined by a formula 
\begin{gather}
\label{2021-06-22-11}
\xi_i=\varphi_i(p_1,p_2,\ldots,p_\sigma)\qquad (i\in{\boldsymbol I}_\tau).
\end{gather}
Then the line element on ${\mathfrak{R}}$ defined by 
\begin{gather}
\label{2021-06-22-12}
ds^2:=d\xi_1^{\;2}+d\xi_2^{\;2}+\cdots+d\xi_\tau^{\;2}=\sum_{\lambda,\mu\in{\boldsymbol I}_\sigma} B_{\lambda,\mu}\,dp_\lambda\,dp_\mu 
\end{gather}
is invariant under an orthogonal transformation $R\in O(\tau)$ leaving stable the submanifold ${\mathfrak{R}}$. 
Further define a volume element on ${\mathfrak{R}}$ using square-root of the discriminant of the quadratic form (\ref{2021-06-22-12}) as 
\begin{gather}
\label{2021-06-22-14}
 \sqrt{|B_{\lambda,\mu}|}\,dp_1\,dp_2\cdots dp_\sigma\,.
\end{gather}
Then it is invariant under $R$ too.

Hurwitz applied this generality to the case of $G=SO(n)\subset {\boldsymbol R}^{n^2}$ with the defining equation (\ref{2021-06-21-2}), where $\tau=n^2,\;\sigma=n(n-1)/2$. 
For $1\le \alpha<n$, denote by $E_\alpha(\varphi)$ the 2-dimensional rotation of angle $\varphi$ in the space of $(x_\alpha,\,x_{\alpha+1})$ given as 
\begin{gather}
\label{2021-06-22-15}
\qquad
\left\{
 \begin{array}{lll}
x_\alpha&=& \;\;\cos\varphi\;x'_\alpha + \sin\varphi\;x'_{\alpha+1}, 
\\
x_{\alpha+1}&=& -\sin \varphi\;x'_\alpha+\cos\varphi\;x'_{\alpha+1},
\\
x_\beta&=& \quad x'_\beta\hspace{18ex} (\beta\ne \alpha,\;\alpha+1).
\end{array}
\right.
\end{gather}

Introduce $n(n-1)/2$ angles as 
\begin{gather}
\varphi_{0,1}\,;\;\varphi_{0,2},\;\varphi_{1,2}\,;\; 
\varphi_{0,3}\,,\;\varphi_{1,3}\,,\;\varphi_{2,3}\,;\;\ldots\,;\;
\varphi_{0,n-1}\,,\;\varphi_{1,n-1}\,,\;\cdots,\;\varphi_{n-2,n-1}
\end{gather}
and take orthogonal transformations as 
\begin{gather}
\label{2021*06-22-16}
\left\{
 \begin{array}{lll}
E_1&=& E_{n-1}(\varphi_{01}), 
\\
E_2&=& E_{n-2}(\varphi_{12})E_{n-1}(\varphi_{02}),
\\
E_3&=& E_{n-3}(\varphi_{23})E_{n-2}(\varphi_{13})E_{n-1}(\varphi_{03}),
\\
\cdots&=& \qquad \cdots\cdots\cdots
\\
E_{n-1}&=& E_1(\varphi_{n-2,n-1})E_2(\varphi_{n-3,n-1})\cdots E_{n-1}(\varphi_{0,n-1}),
\end{array}
\right.
\end{gather}
and put 
\begin{gather}
\label{2021-06-22-17}
R=E_1E_2\cdots E_{n-1}.
\end{gather}
When the parameter $\varphi_{rs}\;(0\le r<s<n)$ runs over 
\begin{gather}
\label{2021-06-23-18}
0\leqq \varphi_{0s}<2\pi,\quad 0\leqq \varphi_{rs}<\pi\;\;(1\le r<s<n),
\end{gather}
the element $R=(r_{\alpha\beta})_{\alpha,\beta\in{\boldsymbol I}_n}$ runs over $G$ once. 

Put\; $ds^2=\underset{\alpha,\beta}{\sum}(dr_{\alpha\beta})^2$\; as the line element on G.

The main results in this case are the following: 
\begin{gather}
\label{2021-06-23-20}
dR=2^{\frac{n(n-1)}{4}}\prod_{0\leqslant r< s<n}(\sin\varphi_{rs})^r\,d\varphi_{rs}\,,\\
\label{2021-06-23-21}
M=\int_GdR=\frac{2^{\frac{(n-1)(n+4)}{4}}\cdot \pi^{\frac{n(n+1)}{4}}}
{\Gamma(1/2)\Gamma(2/2)\Gamma(3/2)\cdots\Gamma(n/2)}\,. 
\end{gather}

These formulas are very interesting, but they were not cited, except by Schur, in modern papers which I read until now.

\subsection{Another form of global coordinates on $SO(n)$}

The above explicit form of global coordinates is not so familiar with me, and taking into account the later result by Schur, I would like to give another form of it well visible from a geometric point of view. For $i,j,\in{\boldsymbol I}_n,\,i\ne j$, take two-dimensional rotation $r_{ij}(\phi)$, in the space spanned by $(x_i,x_j)$, given as 
\begin{gather}
\nonumber
r_{ij}(\phi)\,:\qquad \begin{pmatrix}x'_i \\ x'_j\end{pmatrix}
\begin{pmatrix}\cos\phi&-\sin\phi \\ \sin\phi & \cos\phi\end{pmatrix}
\begin{pmatrix}x_i \\ x_j\end{pmatrix}.
\end{gather}
Put $G_n=SO(n)$ and consider a descending series of subgroups $G_n\subset G_{n-1}\subset \cdots G_3\subset G_2$\, and the corresponding series of quotient spaces $S^{n-1}=G_n/G_{n-1}\,,\; S^{n-2}=G_{n-1}/G_{n-2}\,,\\ 
\cdots\,,\; S^2=G_3/G_2$\,.  
Take a rotation $r_{j-1}\in G_j\;(2\le j\le n)$ given below as
\begin{eqnarray}
\nonumber
\begin{array}{c}
r_{n-1}=r_{12}(\phi_{1,n-1})r_{23}(\phi_{2,n-1})\cdots\cdots r_{n-1,n}(\phi_{n-1,n-1}),
\\[.3ex]
r_{n-2}=r_{12}(\phi_{1,n-2})r_{23}(\phi_{2,n-2})\cdots r_{n-2,n-1}(\phi_{n-2,n-2}),
\\[.3ex]
\cdots\;=\qquad\cdots\cdots\cdots\cdots\qquad, 
\\[.3ex]
\cdots\;=\qquad\;\cdots\cdots\cdots\;\qquad, 
\\[.3ex]
r_2=r_{12}(\phi_{1,2})r_{23}(\phi_{2,2}),
\\[.3ex]
r_1=r_{12}(\phi_{1,1}).
\end{array}
\end{eqnarray}
Then we have the following 
\vskip1.2em

{\bf Theorem 1.2.}\; Put a general element $u$ in $G_n:=SO(n)$ as 
\begin{eqnarray}
\nonumber
u=r_{n-1}r_{n-2}\cdots r_2r_1.
\end{eqnarray}
For a fixed $j\;(1\le j< n)$, when the parameter $\phi_{1,j},\phi_{2,j},\cdots,\phi_{j,j}$ of $r_j$ runs over 
\begin{eqnarray}
\nonumber
0\le \phi_{i,j}\le\pi\;(1\le i<j),\quad 0\le \phi_{j,j}<2\pi,
\end{eqnarray}
it gives a spherical coordinates on $j$-dimensional sphere 
 $S^j\cong G^{j+1}/G^j$. From the total point of view over $1\le j<n$, the radius of sphere decreases smaller and smaller depending on $j$, and when $(\phi_{ij})_{1\leqslant i\leqslant j<n}$ moves over the above region, the total image of $u$ above covers the whole of $G_n=SO(n)$. 
\vskip1.2em

{\it Proof.}\; Let $({\boldsymbol e}_1,{\boldsymbol e}_2,\cdots,{\boldsymbol e}_n)$ be the standard orthonormal frame in a Euclidean space $E^n$ around the origin and put $({\boldsymbol e}'_1,{\boldsymbol e}'_2,\cdots,{\boldsymbol e}'_n),\,{\boldsymbol e}'_i:=u{\boldsymbol e}_i$. Then $u$ is determined by this frame. 
As is known rather well, the unit vector ${\boldsymbol e}'_n$ is obtained as $r_{n-1}^{\;\;\;-1}{\boldsymbol e}'_n={\boldsymbol e}_n$ with an appropriate parameter $(\phi_{i,n-1})_{1\leqslant i\leqslant n-1}$. 
Take a new frame $({\boldsymbol e}''_1,{\boldsymbol e}''_2,\ldots,{\boldsymbol e}''_n),\, {\boldsymbol e}''_i=r_{n-1}^{\;\;\;-1}{\boldsymbol e}'_i,$, then we have ${\boldsymbol e}''_n={\boldsymbol e}_n$. Here, replacing $S^{n-1}\subset E^n$ by $S^{n-2}\subset E^{n-1}$, we repeat the same discussion as above, that is, choose $r_{n-2}$ so that $(r_{n-2})^{-1}{\boldsymbol e}''_{n-1}={\boldsymbol e}_{n-1}$. 
Again we choose a new frame ${\boldsymbol e}^{(3)}_i:=(r_{n-2})^{-1}{\boldsymbol e}''_i\;(i\in{\boldsymbol I}_n)$ in such a way that ${\boldsymbol e}^{(3)}_i={\boldsymbol e}_i\;(i= n-1,n)$. 
This time, we discuss on $S^{n-3}\subset E^{n-2}$. Thus, repeating this process, we arrive at the last stage such that $(r^{\;-1}_1r^{\;-1}_2\cdots r^{\;\;\;-1}_{n-1}){\boldsymbol e}'_i={\boldsymbol e}_i\;(i\in{\boldsymbol I}_n)$. This gives $u=r_{n-1}\cdots r_2r_1$. 
\hfill 
$\Box$

\subsection{On \lq\,\S 3\; Global coordinates on $SU(n)$ and formula of invariant measure\,\rq}

Now return to the original \cite{Hur1897}. 
Denote by $a^0$ the conjugate of a complex number $a$. 
Then an element 
$T=(c_{\alpha\beta})$ of the group $G=SU(n)$ is defined 
as
\begin{gather}
\label{2021-06-21-31}
\sum_\gamma c_{\gamma\alpha}c^0_{\gamma\beta}=\delta_{\alpha\beta}\;;\quad |c_{\alpha\beta}|=1.
\end{gather}

Introduce a two-dimensional unitary transformation in $(x_\alpha,x_{\alpha+1})$-space as 
\begin{gather}
\label{2021-06-21-32}
\left\{
\begin{array}{lll}
x_\alpha&=&\hspace{2ex}ax'_\alpha+bx'_{\alpha+1}\,,
\\
 x_{\alpha+1}&=&-b^0x'_\alpha+a^0x'_{\alpha+1}\,,
\\
x_\beta&=&\quad x'_\beta\,, 
\end{array}
\right.
\\
\label{2021-06-21-33}
a=\cos \varphi\,e^{\psi i},\;\; b=\sin \varphi\,e^{\chi i},\;\;a^0=\cos\varphi\,e^{-\psi i},\;\; b^0=\sin \varphi\,e^{-\chi i},
\end{gather}
which gives all transformation in $SU(2)$. 
Denote it by $E_\alpha(\varphi,\psi,\chi)$, and put 
\begin{gather}
\label{2021-06-21-34}
\left\{
\begin{array}{lll}
E_1&=&E_{n-1}(\varphi_{01},\psi_{01},\chi_1)\,,
\\
E_2&=&E_{n-2}(\varphi_{12},\psi_{12},0)E_{n-1}(\varphi_{02},\psi_{02},\chi_2)\,,
\\
E_3&=&E_{n-3}(\varphi_{23},\psi_{23},0)E_{n-2}(\varphi_{13},\psi_{13},0)E_{n-1}(\varphi_{03},\psi_{03},\chi_3)\,, 
\\
\cdots&=&\qquad\cdots\cdots\,,
\\
E_{n-1}&=&E_1(\varphi_{n-2,n-1},\psi_{n-2,n-1},0)E_2(\varphi_{n-3,n-1},\psi_{n-3,n-1},0)\,\cdots
\\
&&\cdots\cdots\;E_{n-2}(\varphi_{1,n-1},\psi_{1,n-1},0)E_{n-1}(\varphi_{0,n-1},\psi_{0,n-1},\chi_{n-1})\,,
\end{array}
\right.
\\[1ex]
\label{2021-06-21-35}
T=E_1E_2E_3\cdots E_{n-1}\,,
\end{gather}
where the parameters run over the region 
\begin{gather}
\label{2021-06-21-36}
0\leqq \varphi_{rs}<\dfrac{\pi}{2}\,,\quad 0\leqq  \psi_{rs}<2\pi,\quad 0\leqq \chi_s<2\pi.
\end{gather}
Then, by simple calculation, we have 
\begin{gather}
\nonumber
\begin{array}{lll} 
E_1{\boldsymbol e}'_n &=& 
-\,e^{i\varphi_{0,1}}\sin\psi_{0,1}\,e^{-i\chi_1}{\boldsymbol e}'_{n-1}
+ e^{-i\varphi_{0,1}}\cos\psi_{0,1}\,e^{-i\chi_1}{\boldsymbol e}'_n\,,
\\
E_2{\boldsymbol e}'_n &=& 
e^{i(\varphi_{1,2}+\varphi_{0,2})}\sin\psi_{1,2}\sin\psi_{0,2}\,e^{-i\chi_2}{\boldsymbol e}'_{n-2}
\\
&&
-e^{i(-\varphi_{1,2}+\varphi_{0,2})}\cos\psi_{1,2}\sin\psi_{0,2}\,e^{-i\chi_2}{\boldsymbol e}'_{n-1}
+e^{-i\varphi_{0,2}}\cos\psi_{0,2}\,e^{-i\chi_2}{\boldsymbol e}'_n\,,
\\
\cdots &=& 
\qquad \cdots\cdots\cdots
\\
E_{n-1}{\boldsymbol e}'_n&=&
(-1)^{n-1}e^{i(\varphi_{n-2,n-1}+\cdots+\varphi_{2,n-1}+\varphi_{1,n-1}+\varphi_{0,n-1})}\sin\psi_{n-2,n-1}\sin\psi_{n-3,n-1}\cdots 
\\
&&\quad\,\times \sin\psi_{2,n-1}\sin\psi_{1,n-1}\sin\psi_{0,n-1}\,e^{-i\chi_{n-1}}{\boldsymbol e}'_1\;+
\\
&&
(-1)^{n-2}e^{i(-\varphi_{n-2,n-1}+\cdots+\varphi_{2,n-1}+\varphi_{1,n-1}+\varphi_{0,n-1})}\cos\psi_{n-2,n-1}\sin\psi_{n-3,n-1}\cdots 
\\
&&\quad\,\times \sin\psi_{2,n-1}\sin\psi_{1,n-1}\sin\psi_{0,n-1}\,e^{-i\chi_{n-1}}{\boldsymbol e}'_2
\\
&&\pm\qquad \cdots\cdots\cdots\cdots
\\
&&+\,e^{i(-\varphi_{2,n-1}+\varphi_{1,n-1}+\varphi_{0,n-1})}\cos\psi_{2,n-1}\sin\psi_{1,n-1}\sin\psi_{0,n-1}\,e^{-i\chi_{n-1}}{\boldsymbol e}'_{n-2}
\\
&&-\,e^{i(-\varphi_{1,n-1}+\varphi_{0,n-1})}\cos\psi_{1,n-1}\sin\psi_{0,n-1}\,e^{-i\chi_{n-1}}{\boldsymbol e}'_{n-1}
\\
&&+\,e^{-i\varphi_{0,n-1}}\cos\psi_{0,n-1}\,e^{-i\chi_{n-1}}{\boldsymbol e}'_n\,,
\end{array}
\end{gather}
and see that $E_{n-1}{\boldsymbol e}'_n$ covers once $SU(n)/SU(n-1)\cong B^{n-1}:=\{x\in{\boldsymbol C}^n\,;\,\|x\|=1\}$.  

\vskip.2em

Here for unitary group $SU(n)$, the line element $ds^2$ chosen for $T=(c_{\alpha\beta})$ is 
\begin{gather}
\nonumber
ds^2=\sum_{\alpha,\beta}dc_{\alpha\beta}\,dc^0_{\alpha\beta}\hspace{15ex}
(c^0_{\alpha\beta}:=\overline{c_{\alpha\beta}})
\end{gather}
and the corresponding invariant measure $dT$ in the global coordinates 
$(\varphi_{rs},\,\psi_{rs},\,\chi_s)$ is  
\begin{gather}
\label{2021-06-28-1}
dT=\sqrt{n!}\;2^{\frac{n(n-1)}{2}}\cdot\prod_{r,s}\cos\varphi_{rs}\,
(\sin\varphi_{rs})^{2r+1}d\varphi_{rs}\,d\psi_{rs}\cdot\prod_s d\chi_s\,.
\end{gather}

\subsection{On \lq\,\S 4\; General theory for Lie groups\,\rq}

In this last part pp.86--90 of [Hur1897], I feel that Hurwitz discussed generally for Lie groups, maybe in the direction of something like Maurer-Cartan formula. But I could not follow it up in detail.

\section{Schur's simplification of Representation Theory for finite groups}

In the paper \cite{Sch05}, 
Schur simplified in a very reasonable way the theory of 
characters due to Frobenius.  
As will be explained a little later, Schur's method for finite groups 
can be applied directly to a compact group $G$, if the existence 
of invariant integrals on  $G$ is admitted. 

\subsection{Frobenius' first three papers}  
The world's first foundation of Theory of Representations of Groups 
begun with three papers \cite{Fro1896a}, \cite{Fro1896b} 
and \cite{Fro1897} of Frobenius. 
We survey shortly them as preceding results of Schur.  
\vskip1em

{\bf 2.1.1. Survey of the first paper [Fro1896a].} 

Let ${\mathfrak{H}}$ be a finite group with the identity element $E$, 
$h:=|{\mathfrak{H}}|$\, its order,  
$\alpha$ a conjugacy class,  
$h_\alpha:=|\alpha|$ the order of $\alpha$, and  
$\alpha':=\{A^{-1}\,;\,A\in\alpha\}$. Put, for three conjugacy 
classes $\alpha,\beta,\gamma$, 
\begin{gather}
\label{2021-07-15-1}
h_{\alpha\beta\gamma}:=\sharp\{(A,B,C)\,;\,A\in\alpha,\,B\in\beta,\,
C\in \gamma,\;ABC=E\}
\end{gather}

Let $k$ be the number of 
conjugacy classes, then the set $(\chi_\alpha)$ of $k$ numbers 
is called a {\it Gruppencharakter} 
if it satisfies the following equation: 
\begin{gather}
\label{2021-07-15-2}
\hspace*{12ex}
h_\beta h_\gamma \chi_\beta\chi_\gamma=f
\sum_\alpha h_{\alpha'\beta\gamma}
\chi_\alpha\,,
\hspace{14ex} 
\mbox{\rm ({\bf character equation})} 
\end{gather} 
where $f$ is called as {\it Proportionalit\"atsfactor}. 
Using results in [F51],\footnote{
{\bf [F51]} {\it \"Uber vertauschbare Matrizen}, Sitzungsberichte 
der K\"oniglich Preu\ss ischen Akademie der Wissenschaften 
zu Berlin {\bf 1896}, 
601--614. [F51] denotes the number of paper in Gesammelte 
Abhandlungen.} 
he proved, as the main result in \S 2,  
that there exist exactly $k$ different 
solutions as 
\begin{gather}
\label{2021-07-18-1}
\chi_\alpha=\chi^{(\kappa)}_\alpha,\quad f=f^{(\kappa)}\qquad
(\kappa=0,1,\ldots,k-1)
\end{gather} 
and $k\times k$ type determinant $|\chi^{(\kappa)}_\alpha|\ne 0$.
\vskip1.2em

{\bf Explanation 2.1. In modern languages.}\; 
Taking into account the results in \S 5,  
I can explain as follows. 

Let 
Conj(${\mathfrak{H}}):=\{\alpha\}$ be the set of all conjugacy classes 
of ${\mathfrak{H}}$ and put functions on ${\mathfrak{H}}$ as  
\vskip.6em
\hspace*{18ex}$F^{(\kappa)}(A):=(f^{(\kappa)})^{-1}\chi^{(\kappa)}_\alpha$\qquad for\quad $A\in\alpha\subset{\mathfrak{H}}$ 
\\[1ex]
starting 
from $(\chi^{(\kappa)}_\alpha)$. On the other hand, consider the group 
algebra $K[G]$ of a group $G$ with coefficient in $K={\boldsymbol Z}, {\boldsymbol R}$ or ${\boldsymbol C}$. Then 
we know that the center $K[{\mathfrak{H}}]^o$ is spanned by elements  
 $X_\alpha:=\sum_{A\in \alpha}A\;\;(\alpha\in {\rm Conj}({\mathfrak{H}}))$. The defining 
relations are calculated as 
\begin{eqnarray}
\label{2021-07-26-11}
X_\beta\cdot X_\gamma =\sum_{\alpha\in {\rm Conj}({\mathfrak{H}})}\frac{h_{\alpha'\beta\gamma}}{h_\alpha}\,X_\alpha\,.
\end{eqnarray}
Put $e_\alpha:=h_\alpha^{\,-1}X_\alpha$, then we 
have \;
$h_\beta h_\gamma e_\beta e_\gamma =\sum_\alpha h_{\alpha'\beta\gamma}e_\alpha$. 
So the function $F^{(\kappa)}(A)$ for $A\in {\mathfrak{H}}$ gives a representation 
of the commutative algebra ${\boldsymbol Z}[{\mathfrak{H}}]^o$, because 
\begin{gather}
\nonumber
h_\beta h_\gamma F^{(\kappa)}(B) F^{(\kappa)}(C)=
\sum_\alpha h_{\alpha'\beta\gamma}\,
F^{(\kappa)}(A).
\end{gather}

This explains the birthplace of the character equation, and $F^{(\kappa)}$ 
is normalized as $F^{(\kappa)}(E)=1$ at the identity element $E$. 

For direct relation with characters of 
irreducible linear representations  
of the group ${\mathfrak{H}}$, see the survey below 
of the third paper \cite{Fro1897}. 
\hfill
$\Box$
\vskip1.2em

In \S\S 8--10, all {\it Gruppencharakter}\hspace{.3ex}s are calculated explicitly 
for groups, ${\mathfrak{A}}_4,\;{\mathfrak{S}}_4/{\boldsymbol Z}_2^{\;2},\;{\mathfrak{S}}_4,
\\ {\mathfrak{A}}_5,\;{\mathfrak{S}}_5$\;(in \S 8), and for groups $PSL(2,{\boldsymbol Z}_p)$ with $p$ odd prime 
(in \S\S9--10). 

\vskip1em

{\bf 2.1.2. Survey of the second paper [Fro1896b].} 

This paper consists of Introduction and \S\S 1--12. 
Let ${\mathfrak{H}}$ be a group and $h$ its order. Consider a set of variables $x=(x_P)_{P\in{\mathfrak{H}}}$ with index $P\in {\mathfrak{H}}$, and $h\times h$ type determinant 
$$
\Theta(x):=\det(x_{PQ^{-1}}),
$$
where $P, Q$ run over ${\mathfrak{H}}$. This is called as {\it Gruppendeterminante}, and is the main object to be studied here. 
First $\Theta$ is divisible by a linear form 
$\xi:=\sum_{R\in{\mathfrak{H}}}x_R$\,, as is easily seen.   
In Introduction, main results of the paper 
are summarized, and I explain them below.   

(1) \;The number of prime factors of $\Theta(x)$ is equal to the number $k$ of conjugacy classes of ${\mathfrak{H}}$. Let them be $\Phi^{(\lambda)}(x)\;(1\le \lambda\le k)$ with $f^{(\lambda)}:=\dim \Phi^{(\lambda)}$. Then $\Theta$ is decomposed as 
$$
\Theta=\prod_{1\leqslant \lambda\leqslant k} \big(\Phi^{(\lambda)}\big)^{e^{(\lambda)}}.
$$

(2) \;The exponent $e^{(\lambda)}$ is equal to the dimension $f^{(\lambda)}$.

(3) \;The dimension $f^{(\lambda)}$ divides the order $h$.

(4) \;By a certain linear transformation, $\Phi^{(\lambda)}$ becomes a function 
 of $(f^{(\lambda)})^2$ number of variables. 
 
(5) \;Collecting all such variables over $1\le \lambda \le k$, we obtain just 
$h=\sum_\lambda (f^{(\lambda)})^2$ number of independent variables.

(6) \;Put, for every conjugacy class $\alpha$, $x_A=x_B$ for $A,B\in\alpha$. Then we have just $k$ number of variables. 
For a prime factor $\Phi^{(\lambda)}(x)$, there exists a {\it Charakter}\, $(\chi^{(\lambda)}_\alpha)_{\alpha\in{\rm Conj}({\mathfrak{H}})}$ having the following property: define a function $\chi^{(\lambda)}$ on ${\mathfrak{H}}$ as  
$\chi^{(\lambda)}(A):=\chi^{(\lambda)}_\alpha$ for $A\in\alpha$, and a linear function $\xi^{(\lambda)}$ as  
$$
\hspace*{10ex}
\xi^{(\lambda)}:=\dfrac{1}{f^{(\lambda)}}\,\sum_R \chi^{(\lambda)}(R)x_R\,,
\hspace{10ex}
$$
then \;
$\Phi^{(\lambda)}={\xi^{(\lambda)}}^{f^{(\lambda)}}.$

(7) \;The set of functions $\xi^{(\lambda)},\,1\le \lambda \le k$,\, are mutually linearly independent. 

(8) \;{\it Charakter}\, $(\chi^{(\lambda)}_\alpha)$ determines $\Phi^{(\lambda)}(x)$ for general $x=(x_R)_{R\in{\mathfrak{H}}}$ completely. Studies on $\Theta(x)$ can be reduced to those of $\xi^{(\lambda)}$'s, in particular, when $x_{AB}=x_{BA}\;(A,B\in{\mathfrak{H}})$, \;
$$
\Theta=\prod_\lambda \big(\xi^{(\lambda)}\big)^{{f^{(\lambda)}}^2}.
$$ 

(9) \;The calculation of the above decomposition of degree $h$ can be 
reduced to that of size $k$ determinant for $(\alpha,\beta)\in {\rm Conj}({\mathfrak{H}})\times {\rm Conj}({\mathfrak{H}})$ given as 
\begin{gather*}
\det\left({\sum}_\gamma\; \dfrac{1}{h_\alpha}\,h_{\alpha\beta'\gamma}x_\gamma\right)_{\alpha,\beta}=\prod_\lambda\xi^{(\lambda)}.
\\[2ex]
\mbox{\rm In fact, since}\hspace{10ex}  
h_\beta h_\gamma \chi^{(\lambda)}_\beta \chi^{(\lambda)}_\gamma 
=f^{(\lambda)}\,\sum_\alpha h_{\alpha'\beta\gamma}\chi^{(\lambda)}_\alpha\,,\qquad
\hspace{14ex}
\\
{\rm and}
\hspace{14ex} 
\xi^{(\lambda)}=\dfrac{1}{f^{(\lambda)}}\,\sum_\gamma h_\gamma \chi^{(\lambda)}_\gamma 
x_\gamma\,,
\hspace{17ex}
\\
\mbox{\rm we have}\qquad 
h_\alpha\chi^{(\lambda)}_\alpha \xi^{(\lambda)}=\dfrac{1}{f^{(\lambda)}}\,\sum_\gamma h_\alpha\chi^{(\lambda)}_\alpha h_\gamma \chi^{(\lambda)}_\gamma 
 x_\gamma
\hspace{26ex}
\\
\hspace*{10ex}
=\sum_\gamma \sum_\beta h_{\alpha\gamma\beta'} \chi^{(\lambda)}_\beta x_\gamma
=\sum_\beta\Big(\sum_\gamma h_{\alpha\beta'\gamma}x_\gamma\Big)\chi^{(\lambda)}_\beta, 
\\
\therefore\quad 
\sum_\beta\Big(\sum_\gamma h_{\alpha\beta'\gamma}x_\gamma
-h_{\alpha\beta'}\xi^{(\lambda)}\Big)\chi^{(\lambda)}_\beta
=0,
\hspace{12ex}
\\
\therefore\quad 
\det\Big(\dfrac{1}{h_\alpha}\sum_\gamma h_{\alpha\beta'\gamma}x_\gamma
-\delta_{\alpha\beta'}\,r\Big)_{\alpha,\,\beta'}=0\quad
{\rm for}\;\; r=\xi^{(\lambda)}\;\;(1\le \lambda\le k)\;\; 
\mbox{\rm (eigenvalues)}.
\end{gather*}
\vskip.3em

{\bf Explanation 2.2. Personal impression on [Fro1896a] and [Fro1896b].} 

If we recognize that the matrix determinant $\Theta(x)$ 
expresses the left regular representation of a finite 
group $\mathfrak{H}$ on $\ell^2(\mathfrak{H})$, then 
the results from (1) to 
(9) listed above are familiar with us, except the assertion (3), that is, $f$ divides $h$. However, their proofs are all purely algebraic, and very new and refreshing for me, and also very difficult to read.  

So, comparing to my personal self-education process in my student age, 
this kind of composition of the theory of group representations is rather astonishing and I never knows that it is possible, I mean, 
the order of introducing: firstly character equation and secondly matrix determinant, and linear representation of $\mathfrak{H}$ is hidden in the backyard.

\vskip1.2em

{\bf 2.1.3. Survey of the third paper [Fro1897].} 

In \S 1, {\it Charakter} of a Group ${\mathfrak{H}}$ is redefined as a function $\chi(R)$ on ${\mathfrak{H}}$ by a system of equations as follows: with $h=|{\mathfrak{H}}|$, 
\begin{gather}
\label{2021-07-26-1}
\chi(E)=f,
\\
\chi(AB)=\chi(BA)\qquad(A, B\in{\mathfrak{H}}), 
\\
h\,\chi(A)\chi(B)=f\,\sum_R\chi(AR^{-1}BR),
\\
h=\sum_R \chi(R)\chi(R^{-1}). 
\end{gather}

In \S 2, {\it Darstellung durch lineare Substitutionen}\, of a group ${\mathfrak{H}}$ is defined as a correspondence of $A\in {\mathfrak{H}}$ to a regular matrix $(A)$ satisfying $(A)(B)=(AB)$ for $A,B\in{\mathfrak{H}}$. Introducing variables $x_R\;(R\in{\mathfrak{H}})$, Frobenius associates a determinant $F(x):=|\sum (R)x_R|$ and studies its factorization into prime factors, and assert that 

\lq\lq\,{\it every prime factor appears in {\it Gruppendeterminante}\, $\Theta(x)$}.\rq\rq 

A Darstellung $(R)$ is called {\it prime}\, if the polynomial $F(x)$ is prime. 

\vskip.5em 

In \S 4, from a representation $(R)$, its (trace) character $\chi(R)$ is defined as 
$\chi(R):={\rm tr}(R)$. 

\vskip.5em

In \S 6, the relation of characters $\chi$ with linear representations $(R)$ is clarified. For any conjugacy class $\alpha$ of ${\mathfrak{H}}$, the operator $J_\alpha:=\sum_{A\in \alpha}(A)$ commutes with any $(R),\,R\in{\mathfrak{H}}$. 

So, if the representation is prime (or irreducible), $J_\alpha$ must be a scalar multiple of the identity operator $(E)$, and so a number $\chi_\alpha$ can be defined as  
\begin{gather}
\label{2021-07-26-12}
f\,J_\alpha=h_\alpha\chi_\alpha\cdot(E)\qquad {\rm or}\qquad 
J_\alpha=\dfrac{1}{f}\,h_\alpha\chi_\alpha\cdot(E).
\end{gather}  

On the other hand, the relation (\ref{2021-07-26-11}) gives us 
\begin{gather}
\label{2021-07-26-13}
J_\beta\cdot J_\gamma=\sum_\alpha \frac{h_{\alpha'\beta\gamma}}{h_\alpha}\,J_\gamma\;\quad
\mbox{\rm and so}\;\quad 
h_\beta\chi_\beta\cdot h_\gamma\chi_\gamma 
=f\,\sum_\alpha \,h_{\alpha'\beta\gamma}\,h_\alpha\chi_\alpha\,.
\qquad
\end{gather} 

This last equation is just equal to the character equation 
in (\ref{2021-07-15-2}). 
\vskip1.2em 

{\bf My personal impression on the third paper [Fro1897]}. 

In this paper, the fundamental notions of linear representations 
and their (trace) characters are introduced. The relations with group 
algebra $K[\mathfrak{H}]$ and its center $K[\mathfrak{H}]^o$ and also with 
intertwining operators are studied, always in {\it algebraic way}. 
Important parts are closely related to the results in the previous 
two papers. So, in total, I feel, 
the side of treatments from the stand points of 
functional analytic ways is weak, and is still basically under way of 
 developing.

\subsection{Schur's simplification of Character Theory for finite groups}

As explained above, the above three papers of Frobenius are difficult to read 
and not easy to approach. Schur, a pupil of Frobenius, wrote the paper 
\cite{Sch05}, wishing to dissolve them, apart from efforts of Burnside.  
Schur changed the order of introducing notions, first {\it linear representations}, not necessary irreducible, and then their {\it characters}\, and {\it irreducible decompositions}\, and so on, as is in modern text books. Thus he aimed to facilitate and transparentize the theory. This method is better, not only for finite groups but also for general groups. 

Schur himself wrote at the beginning of the paper as  
\begin{quotation}
\noindent 
{\Large D}ie vorligende Arbeit enth\"alt eine durchaus elementare Einh\"urung in die 
Hrn.\\ 
{\sc Frobenius} begr\"undete Theorie der Gruppencharaktere, die auch als die Lehre von Darstellung der endlichen Gruppen durch lineare homogene Substitutionen bezeichnet werden kann.
\end{quotation} 

The paper consists of Preface and Sections 1 to 6 and well-organized. It contains 17 Assertions written in {\it italic}\, and numbered from I to XVII, and Important Equalities numbered as (I.)\,$\sim$\,(III.), (III$'$.), 
(VI.)\,$\sim$\,(VI), (VI$'$.), (VII.)\,$\sim$\,(XIV.)

We notice three important results of general character. 

First one is the so-called Schur's lemma on intertwining operators 
between two irreducible representations. This corresponds 
to Assertions I and II. 

Second one is the so-called {\it orthogonality 
relation}\, for matrix elements of irreducible representations. This corresponds Equalities (I.) and (II.). 

Third one is also {\it orthogonality 
relation}\, among irreducible characters.\!\!
\footnote{trace characters of irreducible linear representations.}  
This corresponds to Equalities (IX.), (X.) and (XII.). 

Let us explain a little more in detail. Schur used repeatedly the averaging of 
a function $f$ on a finite group $G$. Despite of losing the good flavor of the classics, I use intentionally modern notation of 
integration in such a way that
\begin{gather}
\label{2021-07-28-1}
\int_G f(g)\,dg:=\frac{1}{|G|}\sum_{g\in G}f(g).
\end{gather}
Here $dg$ denotes the normalized invariant measure on $G$. 
In the space of functions ${\boldsymbol C}(G)$ on $G$, introduce  
inner product as 
\begin{gather}
\label{2021-07-28-2}
\langle f_1, f_2\rangle :=\int_G f_1(g)\,\overline{f_2(g)}\,dg\,,
\end{gather}
then we get a Hilbert space $L^2(G)$. 
In the main part of the paper, Schur proved, in principle, the following assertions (which I write in modern languages in a form of theorem, for the sake of brevity): 
\vskip1em

{\bf Theorem 2.3.}\; 
(i)\; {\it 
Let $\pi$ be a linear representation of $G$ on a vector space  $V(\pi)$. Then an inner product can be introduced in $V(\pi)$ in such a manner that $\pi(g),\,g\in G$, are all unitary, that is, $\pi$ becomes a unitary representation.
} 

(ii)\; {\it 
In $V(\pi)$, take a complete orthonormal system and, with respect to it, express $\pi(g)$ by a unitary matrix $\big(t^\pi_{\alpha\beta}\big)_{1\leqslant\alpha,\beta \leqslant d_\pi}$, with $d_\pi=\dim \pi$ (\,$t^\pi_{\alpha\beta}$ are matrix elements). 
Suppose $\pi$ is irreducible, and take another 
irreducible unitary representation $\rho$, not equivalent to $\pi$. Then
\begin{gather}
\label{2021-07-28-3}
\langle t^\pi_{\alpha\beta}, t^\pi_{\gamma\delta}\rangle =\dfrac{1}{d_\pi}\,\delta_{\alpha,\beta}\,\delta_{\gamma,\delta}\qquad(1\le \alpha,\beta,\gamma,\delta\le d_\pi),
\\
\label{2021-07-28-4}
\hspace*{.6ex}
\langle t^\pi_{\alpha\beta}, t^\rho_{\gamma\delta}\rangle =0\qquad (1\le \alpha,\beta\le d_\pi,\;1\le \gamma,\beta\le d_\rho).
\end{gather} 
} 
(iii)\; {\it The character of linear representation $\pi$ is defined as 
\;$\chi_\pi(g):={\rm tr}\big(\pi(g)\big)\;\,(g\in G)$, which is  
invariant under inner automorphisms. Suppose $\pi$ 
is irreducible and $\rho$ is another irreducible representation not equivalent 
to $\pi$. Then 
\begin{gather}
\label{2021-07-28-5}
\|\chi_\pi\|^2=\langle \chi_\pi,\chi_\pi\rangle =1,
\qquad \langle \chi_\pi,\chi_\rho\rangle =0\quad
{\rm or}\quad \chi_\pi\perp\chi_\rho\,,
\end{gather}
where $\|f\|$ denotes the norm of $f$ in $L^2(G)$.
}
\vskip1em

I found that Schur's 
methods themselves for proving these assertions, can be directly 
applied to any group $G$ if it has a {\it finite invariant measure}\, on it. 
For instance, take the equalities (\ref{2021-07-28-3}) and 
(\ref{2021-07-28-4}) above, which correspond respectively Equalities (I.) and 
(II.) in \cite{Sch05}. I will rewrite Schur's proof on two pages of 
\cite{Sch05} 
by using integral symbol $\int_Gf(g)dg$ in place of $\sum_{g\in G}f(g)$. 
Then it looks like as follows: 
\begin{quotation}
\noindent
Let $\pi$ be an irreducible matrix 
 representation of dimension $f\;(=d_\pi)$.
\\  
Take an arbitrary matrix $U=(u_{\alpha\beta})$ of degree $f$, and put 
\begin{gather*}
{\rm (3.)}\hspace{15ex}
V=\int_G\pi(g^{-1})U\pi(g)\,dg\,.
\hspace{33ex}
\end{gather*}
Then, for any $g_0\in G$, \;$\pi(g_0^{\,-1})V\pi(g_0) 
=\int_G\pi\big(gg_0)^{-1}\big)U\pi(gg_0)\,dg$. 
Since the measure $dg$ is right-invariant, this integral is equal to \;$\int_G
\pi(g^{-1})U\pi(g)\,dg=V$, whence  
\;$\pi(g_0^{\,-1})V\pi(g_0)=V$\; for any $g_0\in G.$ Therefore, 
by Schur's lemma, 
$$
V=vE_f
$$ 
with scalar $v$ and $f$-dimensional identity 
matrix $E_f$. Using the arbitrariness of $(u_{\alpha\beta})$, 
we can obtain the equality (\ref{2021-07-28-3}) easily.
\end{quotation}  

I emphasize that this method of proof can be applied directly 
to any compact group $G$, 
if once the existence of invariant measure on it with finite volume 
is known. 
Similar comment is valid for the orthogonality (\ref{2021-07-28-4}) and (II.). 

When I recognized this fact, the generality of Schur's method, first time, I was very much impressed by his originality 
despite of his comment \lq\lq durchaus elementare Einf\"urung\rq\rq\,  
in the top line of Introduction.

\section{Schur's work and communications with Weyl}

\subsection{Schur's classification of irreducible 
representations of $n$-th 
rotation groups as application of invariant integral}

After 27 years from Hurwitz's paper \cite{Hur1897}, there appeared 
as its application Schur's papers \cite{Sch24a} and \cite{Sch24b}.

\vskip1em

{\bf 3.1.1\; On the first paper \cite{Sch24a}.} 
\vskip.4em 

The contents are, firstly  
 a general theory of the method of applying invariant integrals on groups 
 to the theory of invariant functions or forms, and secondly a 
general explanation for the  
 method of classification of irreducible representations 
of $n$-th orthogonal groups $\mathfrak{D}\;(:=SO(n))$. 
Schur quotes Hurwitz's result and 
ameliorates it for a global 
coordinates on $\mathfrak{D}$ and for 
explicit expression of invariant measures. 

The second part is carried out by proving orthogonality 
relations of matrix elements, an 
 extension of Theorem 2.3\,(ii) for finite groups to 
the rotation group $\mathfrak{D}$.  
Then Schur determines irreducible characters explicitly by 
using the orthogonality relations among them, which is also an 
extension to $\mathfrak{D}$ of Theorem 2.3\,(iii) for finite groups.   

To explain a little more on the contents, I first made a table of contents 
as below\,: 
\vskip.5em

\quad {\bf 
 Erster Teil. Projektive Invarianten.}

\S 1. Ein Hilfssatz \"uber unit\"are Substitutionen. 

\S 2. Der Integrationsproze\ss\ zur Erzeugung projektiver Invarianten. 

\S 3. Beziehungen zum $\Omega$-Proze\ss.

\quad {\bf 
 Zweiter Teil. Die Homomorphismen\/\footnote{
Here \lq\lq\,Homomorphism\,\rq\rq\/ means linear representation.}
 der Drehungsgruppe und das 
\\ \hspace*{16.2ex} 
Abz\"arungsproblem f\"ur Orthogonalinvarianten. }

\S 4. {\sc Hurwitz}sche Integralkalk\"ul. 

\S 5. Einige Eigenschaften der Homomorphismen der Gruppe ${\mathfrak{D}}.$

\S 6. Die Grundrelationen f\"ur die enfachen Charakteristiken.\footnote{
Here \lq\lq\,Charakteristik\,\rq\rq\/ means (trace) character.}

\S 7. Das Abz\"ahlungsproblem f\"ur Orthogonalinvarianten. 

\S 8. Die F\"alle $n=2$ und $n=3$. 

\S 9. Beliebige orthogonale Transformtionen. 

\vskip.6em 

In \S 9, there appear discussions for the full orthogonal groups 
${\mathfrak{D}}'\;(:=O(n))$ too. 

In addition to this, I will quote the last paragraph of his Introduction: 
\begin{quotation}
\hspace*{.2em}
Im Falle der Orthogonalinvarianten sheint der 
{\sc Hurwitz}sche Integrations-\\
proze\ss\/ 
keine \"ahnliche Umgestaltung zuzulassen. 
Ich zeige aber, da\ss\/ gerade im diesen Falle der von 
{\sc Hurwitz} entwickelte Kalk\"ul noch andere wichtige 
Anwendungen gestattet. 
Insbesondere liefelt er eine elegante L\"osung des 
\lq\lq\,Abz\"ahlungs-\\
problems\,\rq\rq, n\"amlich der Aufgabe, 
die genau Anzahl der zu $f(a,x)$ geh\"orenden linear 
unhabh\"angigen Orthogonalinvarianten von vorgegebenem 
Grade $r$ zu bestimmen. Man gelangt zu dieser L\"osung auf dem 
im Falle einer endlichen Gruppe ${\mathfrak{G}}$ vom {\sc Molien} eingeschlagen Wege, 
indem man das an und f\"ur sich wichtige Studium der 
mit der \lq\lq\,Drehungsgruppe\,\rq\rq\/ ${\mathfrak{D}}$ homomorphen 
Gruppen linearer homogener Substitutionen weiterverfolgt 
und eine Theorie entwickelt, die weitgehende 
Analogien mit der sch\"onen {\sc Frobenius}schen 
Theorie der Gruppencharaktere aufweist. 
\end{quotation} 

I understand, in short, that the way of calculation of Hurwitz 
plays an important role in \lq\,Counting-up Problem,\rq\/ for $f(a,x)$ 
for the orthogonal group 
$\mathfrak{D}$, 
and also in \lq\,Classification Problem\,\rq\/ of 
irreducible linear representations for 
$\mathfrak{D}$ through an extension of Frobenius' character theory 
from the case of finite groups to that of compact Lie groups 
$\mathfrak{D}$.

\vskip1.2em

{\bf 3.1.2\; On the second paper \cite{Sch24b}.} 
\vskip.4em 

In this paper Schur actually succeeded to classify irreducible 
linear representations of the rotation group ${\mathfrak{D}}=SO(n)$ and the full 
orthogonal group ${\mathfrak{D}}'=O(n)$. To arrive to 
these results, the theory of characters plays a decisive role. 
The paper contains Satz I to Satz X, printed in {\it italic}. 
To explain the contents, I made table of contents shown below: 
\vskip.5em

\hspace*{1ex}
After two and a half pages of Introduction, 

\S 1. Allgemeine Vorbemerkungen. 

\S 2. Die F\"alle $n=2$ und $n=3$. 

\S 3. Eine Hilfsbetrachtung. 

\S 4. Die einfachen Charakteristiken der Gruppe ${\mathfrak{D}}'$. 

\S 5. Fortsetzung und Schulu\ss\/ des Beweises.

\S 6. Folgerungen aus dem Satze IV.  

\vskip.5em

Quoting the third paragraph of Introduction, we 
continue to explain the contents of the paper: 
\begin{quotation}
Mein Hauptergebnis lautet: Alle stetigen Homomorphismen%
\footnote{
In these ages, Schur put an equal footing for a linear 
representation (say $\pi$) of $G$ and its image $\pi(G)$, and the word 
\lq\lq\,Homomorphism der Gruppe ${\mathfrak{D}}$\,\rq\rq\/ means linear representation 
of ${\mathfrak{D}}$. He does not use any symbol (such as $\pi$) to denote an irreducible 
linear representations, but, in place of such symbols, he uses images of it, 
for instance for $n$ odd, to denote irreducible representations associated 
with ${\boldsymbol \alpha}$, he uses the symbol \lq\lq\,eine mit 
${\mathfrak{D}}$ homomorphe 
irreduzible Substitutionsgruppe ${\mathfrak{G}}_{\alpha_1,\alpha_2,\ldots,\alpha_\nu}$ in 
$N^{(n)}_{\alpha_1,\alpha_2,\ldots,\alpha_\nu}$ 
Variabeln\,\rq\rq. 
} 
der Gruppen ${\mathfrak{D}}$ und ${\mathfrak{D}}'$ lassen sich allein unter 
Benutzung ganzer rationaler Funktionen herstellen. 
Um eine \"Ubersicht \"uber die Gesamtheit der irreduziblen Gruppen linearer 
homogener Substitutionen zu gewinnen, die der Gruppe ${\mathfrak{D}}$ bzw. ${\mathfrak{D}}'$ homomorph sind, verfahre man folgenderma\ss\/en. 
Bedeutet $\nu$ die Zahl $\Big[\dfrac{n}{2}\Big]$, 
so bilde man f\"ur jedes System von $\nu$ nicht negativen ganzen Zahlen 
\begin{gather*}
(1)
\hspace*{23ex}
\alpha_1 \geqq \alpha_2\geqq \cdots\geqq \alpha_\nu
\hspace{32ex}
\end{gather*}
$\cdots\cdots\cdots$ (I continue the explanation in my words) $\cdots\cdots\cdots$
\end{quotation} 

Schur gave explicitly dimension formulas for irreducible linear 
representations of ${\mathfrak{D}}$ and ${\mathfrak{D}}'$ corresponding to 
${\boldsymbol \alpha}:=(\alpha_1,\alpha_2,\ldots,\alpha_\nu)$, clearly 
depending on the parity of $n$, with commentary explanations in detail. 
He continued as 
\begin{quotation}
Bei dem Beweis der hier erw\"ahnten S\"atze benutze ich die auf dem 
{\sc Hurwitz}schen Integralkalk\"ul beruhenden analytischen Methoden, 
die ich in A.I (=[Sch24a]) entwickelt habe. Es handelt sich in erster Linie darum, 
gewisse Integrale mit wohlbestimmten Integranden genau zu berechnen. 
 $\cdots\cdots\cdots\cdots\,$ 
\end{quotation} 

In the next paragraph, Schur commented two papers of \'E. Cartan\footnote{
\'E.\;Cartan, 
{\it Les groupes projectifs qui ne laissent invariants aucune multiplicit\'e 
plane}, 
Bull. Soc. Math. France, {\bf 41}(1913), 53--96, and \;---, 
{\it Les groupes projectifs continues r\'eels qui ne laissent 
invariants aucune multiplicit\'e plane}, 
J. de Math., S\'er. 6, {\bf 10}(1914), 149--186.
}  
 and thanked H.\;Weyl for a notice on these Cartan's works. 

\vskip1.2em 

{\bf A note from my point of view.}\; 
The group ${\mathfrak{D}}=SO(n),\,n\ge 3,$ has universal covering group Spin$(n)$, and so ${\mathfrak{D}}$ has two-valued representations, called {\it spin} and their highest 
weights ${\boldsymbol \alpha}=(\alpha_i)_{i\in{\boldsymbol I}_\nu}$ consist of half-integers (halfs of odd integers) $\alpha_i$'s. In the works of Cartan, linear representations are treated on the level of Lie algebras and there does not encounter with strict differences between {\it spin}\, and {\it non-spin}\, cases. In these days, 
study of covering groups has not yet begun.

\subsection{Weyl's idea inspired by Schur (an advanced announcement)}

{\bf 3.2.1. On the paper [Wey24].} 
\vskip.3em

This paper of Weyl has subtitle 
\lq\lq\,From a letter to Mr. I.\;Schur\,\rq\rq. Actually, 
before this letter, Schur kindly sent to Weyl the drafts of two papers 
[Sch24b] and [Sch24c] before their publication, and Weyl's letter 
is nothing but its reply, from which the contents of this paper are taken.  
Weyl got to the heart of Schur's method, which I can explain as follows:  
\vskip.5em
{\it 
The essential of Schur's method is not global explicit 
expression of invariant measure on the group ${\mathfrak{D}}$ and its application, 
but is 
the integration of invariant functions (such as characters),  
the role of Cartan subgroups, and Weyl groups and the 
so-called Weyl's integration formula, not on the whole of the group\, 
${\mathfrak{D}}$ but 
essentially on Cartan subgroup\, {\rm (cf. {\bf Explanation 4.1} in the next 
section)}. 
}
\vskip.7em

Weyl takes classical groups ${\mathfrak{G}}=SL(n,{\boldsymbol C}),\,{\mathfrak{C}}=Sp(2n,{\boldsymbol C}),$ and  
${\mathfrak{D}}=SO(n,{\boldsymbol C})$, and recognized that it is sufficient to work on 
their compact real forms ${\mathfrak{G}}_u=SU(n),\,{\mathfrak{C}}_u=SpU(2n):=Sp(2n,{\boldsymbol C})\cap U(2n)$,  and 
${\mathfrak{D}}_u=SO(n),$ and on their Cartan subgroups. He explains in detail 
how to classify irreducible unitary representations of these compact 
classical groups by determining 
irreducible characters (the so-called unitarian trick or \lq\lq\,die unit\"are Beschr\"ankung\,\rq\rq). 
 
\vskip1.2em

{\bf 3.2.2. On the third paper [Sch24c].} 
\vskip.3em

In Introduction, there are comments on Theory of Study on rotation 
invariants and Weyl's integration formula on invariant functions in comparison to Hurwitz integral. I made a list of contents as follows: 

\vskip.3em

\qquad \S 1. Einige Hilfsformeln. 

\qquad \S 2. Der vereinfachte Integralkalk\"ul.

\qquad \S 3. Der Abz\"ahlungskalk\"ul f\"ur Orthogonalinvarianten. 

\qquad \S 4. \"Uber die reellen Darstellungen der Gruppe ${\mathfrak{D}}$. 

In addition, with reference to the niceness of the exchanges with Schur and Weyl, one can see about the details of these papers in my article\footnote{
T.\,Hirai, {\it 
On three papers [S51], [S52], [S53] of Schur and a paper [W61] of Weyl
}\, (in Japanese), in 
Proc. of the 27th Symp. on History of Mathematics, 
held October 2016, Report of 
Tsuda University IMCS, {\bf 38}(2017), 50--67.  
https://www2.tsuda.ac.jp$>$math$>$suugakushi$>$sympo27
} 
written in Japanese.

\subsection{Schur's addresses to introduce these articles to Sitzungsberichte}

Academician Schur gave addresses to introduce these articles to 
Sitzungsberichte der Preuss- 
ischen Akademie 
der Wissenschaften zu Berlin {\bf 1924}, and I feel that these addresses give testimony 
about interrelations among Hurwitz, Schur and Weyl. So, I would like to copy them from mimutes of the academy.\footnote{
F.\,Frobenius (1849/10/26 - 1917/08/03),\; 
 A.\,Hurwitz (1853/03/26 - 1919/11/18),
\\
\hspace*{3.6ex}
 I.\,\;Schur \;(1875/01/10 - 1941/01/10),\hspace{4.4ex}
 H.\,\;Weyl \;(1885/11/09 - 1955/12/08).
}  
\vskip1.2em 

{\bf Introductive address of 
[Sch24a] at the meeting on the 1st of January:} 
\begin{quotation}
Der von A.\,Hurwitz angegebene Integralkalk\"ul zur Erzeugung von Invarianten l\"a\ss\/t sich im Falle der projektiven Invarianten durch einen einfacheren Kalk\"ul ersetzen. F\"ur den Fall der Orthogonalinvarianten liefelt die {\sc Hurwitz}sche Methode auch eine L\"osung des Abz\"ahlungsproblems. Dies gelingt, indem man f\"ur die Homomorphismen der Gruppe der reellen orthogonalen Substitutionen eine Theorie entwickelt, die weitgehende Analogien mit der {\sc Frobenius}schen Theorie der Gruppencharaktere aufweist.
\end{quotation}
\vskip.3em

{\bf Introductive address of 
[Sch24b] on the 13th of November:} 
\begin{quotation}
Mit Hilfe der vom Verfasser fr\"uher entwickelten analytischen Methode werden alle irreduziblen Darstellungen der reellen Drehungsgruppe durch lineare homogene Substitutionen n\"aher bestimmt und die zugeh\"origen Variabelnanzahlen genau berechnet. F\"ur die Gruppe aller reellen orthogonalen Substitutionen l\"a\ss\/t sich auch die Gesamtheit aller einfachen Charakteristiken angeben. 
\end{quotation}

\vskip.3em

{\bf Introductive addresses of [Wey24] and 
[Sch24c] on the 11th of December:}
\begin{quotation}
\noindent
{\bf For [Wey24].}

Durch eine Modifikation der von Hrn.\;Schur im Falle der reellen Drehungsgruppe entwickelten Methode gelingt es dem Verfasser, f\"ur alle einfachen und halbeinfachen kontinuierlichen Gruppen das schon von Hrn.\;E.\,Cartan behandelte Darstellungsproblem in abgeschlossenerer Form zu l\"osen.
\end{quotation}

\begin{quotation}
\noindent
{\bf For [Sch24c].} 

F\"ur die reellen Drehungsgruppe wird auf analytischem Wege eine Umgestaltung des {\sc Hurwitz}schen Integralkalk\"uls gewonnen, die, wie Hr.\;Weyl durch eine geometrische Betrachtung gezeigt hat, auch bei allen anderen im Betracht kommenden Gruppen durchgef\"uhrt werden kann. 
\end{quotation}

\section{Weyl's classification of irreducible representations of complex simple Lie groups} 

By the method announced in [Wey24], Weyl succeeded to classify 
(finite-dimensional) irreducible representations of all classical simple 
groups ${\mathfrak{G}}=SL(n,{\boldsymbol C}),\,{\mathfrak{C}}=Sp(2n,{\boldsymbol C}),\, {\mathfrak{D}}=SO(n,{\boldsymbol C})$. Then he applied the same method to complex simple groups of 
exceptional type. The results were published in four parts in 
\cite{Wey25-6}. 

 The most important new invention of Weyl for the step-up from Schur's method 
is the so-called integration formula of Weyl for compact simple Lie groups, which is applicable to determine explicitly irreducible characters. 
Taking its importance into account, 
 I'd like to continue my explanation. 
\vskip1.2em

{\bf Explanation 4.1.}\; 
Let $G$ be a complex simple Lie group and ${\mathfrak{g}}={\rm Lie}(G)$ its Lie algebra. 
Take a compact real form ${\mathfrak{g}}_u$ of ${\mathfrak{g}}$ and let $G_u\subset G$ be the closed 
subgroup of $G$ corresponding to ${\mathfrak{g}}_u$. For instance, for $G={\mathfrak{G}},\,{\mathfrak{C}}$ and ${\mathfrak{D}}$ above, we have $G_u=SU(n),\,SpU(2n)$ and $SO(n)$ respectively. 
To simplify the situation, we take here 
$G={\mathfrak{D}}=SO(n,{\boldsymbol C})$ and $G_u=SO(n)$. 

For $1\le i<j\le n$, let $r_{ij}(\varphi)$ be two-dimensional rotation in $(x_i,x_j)$-space and $A_{ij}$ 
its basis matrix in ${\mathfrak{g}}_u$ 
given as  
\begin{gather}
\label{2021-08-18-1}
r_{ij}(\varphi)=\begin{pmatrix}
\cos\varphi &-\sin\varphi \\ 
\sin \varphi &\cos\varphi
\end{pmatrix}
,\qquad 
A_{ij}=\begin{pmatrix}
0 &-1 \\ 
1 &0
\end{pmatrix}.
\end{gather}
Let ${\mathfrak{h}}_u$ be the Cartan subalgebra of ${\mathfrak{g}}_u$ consisting of elements 
$$
X({\boldsymbol \varphi}):=\varphi_1A_{12}+\varphi_2A_{34}+\cdots+\varphi_\nu A_{2\nu-1,2\nu},\qquad {\boldsymbol \varphi}:=(\varphi_j)_{j\in{\boldsymbol I}_\nu}\,,
$$
and let $H_u$ be the centralizer $Z_{G_u}({\mathfrak{H}}_u)$ of ${\mathfrak{h}}_u$ in $G_u$, then $H_u$ is connected and is a Cartan subgroup of $G_u$ consisting of 
 \;$h({\boldsymbol \varphi}):=
\exp X({\boldsymbol \varphi})=r_{12}(\varphi_1)r_{34}(\varphi_2)\,\cdots\,r_{2\nu-1,2\nu}(\varphi_\nu).
$\, 

Let 
${\mathfrak{h}}$ be the complexification of ${\mathfrak{h}}_u$, and consider root system for 
the pair $({\mathfrak{g}},{\mathfrak{h}})$ and introduce an appropriate order, and let ${\boldsymbol e}_j$ be an element of the dual space ${\mathfrak{h}}^*$ given as $\langle {\boldsymbol e}_j, X({\boldsymbol \varphi})\rangle:=i\varphi_j$, and for ${\boldsymbol \alpha}=(\alpha_1,\alpha_2,\ldots,\alpha_\nu)$, we put ${\boldsymbol e}({\boldsymbol \alpha}):=\sum_{j\in{\boldsymbol I}_\nu}\alpha_j{\boldsymbol e}_j$ and 
\begin{gather}
\label{2021-08-18-21}
\xi_{\boldsymbol \alpha}\big(h({\boldsymbol \varphi})\big):=\exp\langle {\boldsymbol e}({\boldsymbol \alpha}),X({\boldsymbol \varphi})\rangle =\exp\big(\alpha_1\varphi_1+\alpha_2\varphi_2+\cdots+\alpha_\nu\varphi_\nu). 
\end{gather} 

The set $\Sigma^+$ of positive roots and a half of the sum of positive roots ${\boldsymbol e}({\boldsymbol \rho})$ are 
\begin{gather}
\label{2021-08-18-11}
\Sigma^+=\big\{
{\boldsymbol e}_j\pm{\boldsymbol e}_k\; (j<k,\;j,k\in{\boldsymbol I}_\nu),\;\; {\boldsymbol e}_j\;(j\in{\boldsymbol I}_\nu)
\big\},
\quad 
\rho=\Big(\frac{2\nu\!-\!1}{2},\,\frac{2\nu\!-\!3}{2},\,\cdots,\,
\frac{1}{2}\Big),
\\
\nonumber
\hspace*{46.6ex}\mbox{\rm in case\, 
$n=2\nu+1\;(\nu\ge 2),$ or of type $B_\nu$},
\\
\label{2021-08-18-12}
\Sigma^+=\big\{{\boldsymbol e}_j\pm{\boldsymbol e}_k\quad (j<k,\;j,k\in{\boldsymbol I}_\nu)\big\},\qquad 
\rho=(\nu-1,\,\nu-2,\,\cdots,\,1,\,0),
\hspace{1ex}
\\
\nonumber
\hspace*{43ex}\mbox{\rm in case\, 
$n=2\nu\;(\nu\ge 4),$ or of type $D_\nu$}. 
\end{gather}

Introduce the so-called Weyl denominator $D(h),\, h\in H_u,$ as  
\begin{gather}
\label{2021-08-19-1}
D(h):=\xi_\rho(h)\prod_{{\boldsymbol \alpha}\in\Sigma^+}
\big(1-\xi_{-{\boldsymbol \alpha}}(h)\big)
=\prod_{{\boldsymbol \alpha}\in\Sigma^+}
\big(\xi_{{\boldsymbol \alpha}/2}(h)-\xi_{-{\boldsymbol \alpha}/2}(h)\big), 
\end{gather}
where, the most right hand side is a symbolic expression, and in case $n=2\nu+1$,\, $\xi_\rho(h)$ is doubly-valued but becomes single-valued on the universal covering group. 
Moreover, denote by $N_{G_u}(H_u)$ (resp. $Z_{G_u}(H_u)$) the normalizer (resp. centralizer) in $G_u$ of $H_u$. Then $Z_{G_u}(H_u)=H_u$ and the quotient 
$W_{G_u}(H_u):=N_{G_u}(H_u)/Z_{G_u}(H_u)=N_{G_u}(H_u)/H_u$ 
is a finite group called Weyl group, which acts on $H_u$.  

Now take an invariant measure $dg$ on $G_u$, normalized as $\int_{G_u}dg=1$. 
For a continuous function $f$ on $G_u$, since $dg$ is two-sided invariant, we can define an invariant function $f^o$ by  
\begin{gather}
\label{2021-08-19-2}
f^o(g):=\int_{G_u}f(vgv^{-1})\,dv\qquad(g\in G_u).
\end{gather}

On the other hand, any element $g\in G_u$ is conjugate to an $h\in H_u$, and 
so the invariant function $f^o$ is uniquely determined by its restriction on $H_u$. Further, if $g$ is regular,\footnote{By definition, $g$ is regular if 
$sh\ne h\;(\forall s\in W_{G_u}(H_u),\,\ne e)$. 
}  
then $h$ is regular too. 
For a regular element 
$h\in H_u$, there exist $|W_{G_u}(H_u)|$ number of elements $h'=vhv^{-1}\in H_u$. Thus the map 
 $(G_u/H_u)\times H_u\ni (vH_u,h)\mapsto vhv^{-1}\in G_u$ is surjective and (except lower dimensional sets) $|W_{G_u}(H_u)|$-multiple. 

The celebrated integration formula of Weyl for the group $G_u$ is given as 

\vskip1.2em

{\bf Theorem 4.2.}\; For a continuous function $f$ on a compact semisimple Lie group $G_u$, there holds the following integration formula: 
\begin{eqnarray*}
\label{2018-08-09-17}
\int_{G_u}f(g)\,dg=\frac{1}{|W_{G_u}(H_u)|}\int_{H_u}\!\int_{G_u}f(vhv^{-1})\,dv\;|D(h)|^2dh
\\
\nonumber 
=\frac{1}{|W_{G_u}(H_u)|}\int_{H_u}\!f^o(h)\,|D(h)|^2dh, 
\hspace{10.5ex}
\end{eqnarray*}
where $dg, dv$ and $dh$ are the normalized invariant 
measures on $G_u,\,G_u$ and $H_u$.
\hfill
$\Box$
\vskip1.2em

{\bf Personal souvenir.}\; 
Very long ago in my student age, I was studying Weyl's 
 book {\it Classical Groups}\, and found that 
his short proof of Integration formula (corresponding to the above 
Theorem 4.2) was difficult 
to understand 
for me, and felt that it is written by a great Math genius, 
and struggled myself to get my own proof. 
Very later, some years ago, when I read Chern-Chevalley paper 
on \'E.\;Cartan and his work,\footnote{
S.-S.\,Chern and C.\,Chevalley,\, {\it \'E. Cartan and his mathematical work},\, Bull. Amer. Math. Soc., {\bf 58}  
(1952), 217--250.
}
 I found a sentence (p.226), very touched me, as  
\begin{quotation}
$\cdots\cdots\cdots$\; Whereas Weyl's line of attack was, if we may say so, 
brutally global, depending essentially on the method of integration on the whole group, the work of Cartan puts the emphasis on the connection between the local and the global.\; $\cdot\cdots\cdots$
\end{quotation} 
At that moment, I couldn't help smiling since I was feeling the original proof 
of Weyl in the book as quite \lq\,{\it brutally geometric}\,\rq\/ at least for me at that age. Anyhow it 
is a geometric proof as Schur commented \lq\lq\,wie Hr.\;Weyl durch eine geometrische Betrachtung gezeigt hat,\rq\rq\/ in the last sentence of the previous section.  

For more details of \cite{Wey25-6}, see pp.178--182 of 
the paper by Chevalley-Weil.\footnote{
C.\,Chevalley et A.\,Weil,\; 
{\it Hermann Weyl (1885\,--\,1955)}\,, 
L'Enseignement Math\'ematique, {\bf III-3}\;(1957), 
157--187.
}

\section{Theorem of Peter and Weyl}

At the stage when the existence of invariant measures on Lie groups 
was known, Peter and Weyl \cite{PW27} proved the {\it completeness of 
irreducible unitary 
representations for compact Lie groups ${\mathfrak{G}}$}, or irreducible 
decomposition 
of the regular representations on $L^2({\mathfrak{G}})$. I made Table of Contents as  
\vskip.4em

\S 1.\; Grundlagen. Die Orthogonalit\"atsrelationen\hspace{22.3ex} pp.737--741

\S 2.\; {\sc Bessel}sche Ungleichung. Ansatz des Problems \hspace{17.3ex} pp.741--744 

\S 3.\; Konstruktion der h\"ochsten zu einer Gruppenzahl 
\\
\hspace*{7.3ex} geh\"origen Darstellung 
\hspace{43.2ex} pp.744--745 

\S 4.\; Zerf\"allung der gewonnenen Darstellung 
\hspace{25.8ex} 
pp.746--749

\S 5.\; Iteration. Beweis der Vollst\"andigkeitsrelation \hspace{19.7ex}
pp.749--752

\S 6.\; Entwicklungssatz. Approximationssatz. Anwendungen \hspace{10.5ex}
pp.753--755

\vskip.5em

The contents of \S 1 are separated into four Parts as 

1. {\it Gruppe. Volummessung in der Gruppenmannigfaltigkeit.}

2. {\it Darstellung.} \quad 3. {\it Charakteristik}. 
  
4. {\it Die Orthogonalit\"atsrelationen. Jede Darstellung ist einer solchen 
\"aquivalent, deren}

\hspace*{2.8ex}{\it Matrizen $E(s)$ unit\"ar sind.} 
\\
In Part 4, there quoted Schur's results as follows: Take two 
irreducible matrix representations, not mutually equivalent,    
$E(s)=\big(e_{ik}\big)_{i,k\in{\boldsymbol I}_n}$,\;  
$E'(s)=\big(e'_{\iota\kappa}\big)_{\iota,\kappa\in{\boldsymbol I}_{n'}}\;(s\in{\mathfrak{G}})$, then 
\begin{gather}
\label{2021-08-21-1}
\int_{\mathfrak{G}} e_{ik}(s)e'_{\kappa\iota}(s^{-1})\,ds\;=\;0\,,
\\
\label{2021-08-21-2}
\int_{\mathfrak{G}} e_{ik}(s)e_{\kappa\iota}(s^{-1})\,ds\;=\;\frac{1}{n}\,\delta_{i\iota}\delta_{k\kappa}\,,\end{gather}
where $ds$ denotes the normalized invariant measure on ${\mathfrak{G}}$, and $n=\dim E$. 

Note that $ds$ is both left- and right-invariant and also inversion-invariant $d(s^{-1})=ds$. 

The character (\lq Charakteristik\rq) is defined as $\chi(s):={\rm tr}(E(s)),\,\chi'(s):={\rm tr}(E'(s))$. When the representations $E$ and $E'$ are assumed to be unitary, then in the Hilbert space ${\cal H}:=L^2({\mathfrak{G}},ds)$, matrix elements $e_{ik}$ and $e'_{\iota\kappa}$ are mutually orthogonal, and characters $\chi$ and $\chi'$ are unit vectors mutually orthogonal. 
\vskip.8em

The contents of \S 2 is an explanation of the method pursuing analogy with 
Fourier transformation on compact abelian groups, e.g., ${\mathfrak{G}}={\boldsymbol T}^k$. 

For a compact Lie group ${\mathfrak{G}}$, denote by $C({\mathfrak{G}})$ the space of complex-valued 
continuous functions on ${\mathfrak{G}}$ and introduce {\it Multiplikation}\, (convolution product): for $x, y\in C({\mathfrak{G}})$,  
\begin{gather}
\label{2021-08-22-1}
xy(s):=\int x(sr^{-1})y(r)\,dr
\\
\nonumber 
{\rm or}\hspace{26ex} 
 xy(st^{-1}):=\int x(sr^{-1})y(rt^{-1})\,dr\,.
\hspace{28ex}
\end{gather}
An element $x\in C(G)$ is called {\it Gruppenzahl}\; 
(nowadays, $C({\mathfrak{G}})$ is called group algebra). 

\vskip1.2em

{\bf Explanation 5.1.}\; 
In principle, Peter and Weyl consider the right regular representation $R$ on the Hilbert space ${\cal H}$, in which $C({\mathfrak{G}})$ is densely contained. Define for $s\in{\mathfrak{G}}$ the right translation operator as 
$(R(s)f)(u):=f(us)$ for $f\in{\cal H},\,u\in{\mathfrak{G}}$. 
For a continuous kernel function $K(u,v)$ on ${\mathfrak{G}}\times{\mathfrak{G}}$, 
define an integral operator $I_K$ on ${\cal H}$ as follows:  
\begin{gather}
\label{2021-08-21-11}
(I_Kf)(u):=\int_{\mathfrak{G}} K(u,v)f(v)\,dv\qquad{\rm for}\;\;f \in {\cal H}.
\end{gather}
(In this case we can define the {\it trace} of $I_K$ by \;${\rm tr}(I_K):=\int_{\mathfrak{G}} K(u,u)\,du$.)

Suppose $I_K$ commutes with all $R(s)$ (or $I_K$ is an intertwining operator 
for $R$ with itself), that is, \,$R(s)I_K=I_KR(s)$. 
Then \;$K(us,v)=K(u,vs^{-1})$, whence 
\;$K(u,v)=K(uv^{-1}, {\small o})$. Here ${\small o}$ denotes the identity element of ${\mathfrak{G}}$. Put $x(u):=K(u,{\small o})\in C(G)$, $I(x):=I_K$, then  
\begin{gather}
\label{2021-08-21-21}
(I(x)f)(u):=(I_Kf)(u)=\int_{\mathfrak{G}} x(uv^{-1})f(v)\,dv=:(x*f)(u),
\end{gather}
where $x*f$ is the {\it convolution product}, and 
$$
I(xy)=I(x)I(y), \qquad {\rm tr}(I(x))=x({\small o})=:S(x)\quad  {\rm (put)}. 
$$ 

Let $L$ be the left regular representation on ${\cal H}$ defined by 
$(L(s)f)(u):=f(s^{-1}u)$, and assume that $I_K$ commutes with $L(s)$ for 
all $s\in {\mathfrak{G}}$. Then we see similarly as for the above case that 
$K(u,v)=K(v^{-1}u,{\small o})$, and so \;$(I_Kf)(u)=(f*x)(u)$. 

If $I_K$ commutes with both $R$ and $L$, then the kernel function $x$ satisfies \;$x(uv^{-1})=x(v^{-1}u)$\; or \;$x(vuv^{-1})=x(u)\;\;(u,v\in{\mathfrak{G}})$, that is, $x$ is an invariant function or a class function. Of course,  
every character is a class function.

On the other hand, for any matrix representation $E(s)$ of ${\mathfrak{G}}$, we can define a representation of the group algebra $C({\mathfrak{G}})$, with $*$-action $x^*(s):=\overline{x(s^{-1})}$, and its trace as 
\begin{gather}
\label{2021-08-22-11}
E(x):=\int_{\mathfrak{G}} E(s)x(s)\,ds 
\\ 
\nonumber
{\rm and}\hspace{28ex}
 {\rm tr}(E(x))=\int_{\mathfrak{G}}\chi(s)x(s)\,ds. 
\hspace{32ex}
\end{gather}
Moreover define a bounded operator on ${\cal H}$ from the character $\chi$ as 
\begin{gather}
\label{2021-08-22-21}
P_\chi:=\chi({\small o})\int_{\mathfrak{G}} R(s)\chi(s^{-1})\,ds=\dim E\int_{\mathfrak{G}} R(s)\overline{\chi}(s)\,ds, 
\\
\nonumber
{\rm then}\hspace{15ex} 
(P_\chi f)(u)=\dim E\int_{\mathfrak{G}}\chi(us^{-1})f(s)\,ds\quad(f\in{\cal H},\;\dim E=n), 
\hspace{23ex}
\end{gather}
and \;$(P_\chi)^2=P_\chi,\;(P_\chi)^*=P_\chi$. Hence $P_\chi$ is an orthogonal 
projection, and its image $P_\chi{\cal H}\subset {\cal H}$ is spanned by matrix elements $e_{jk}(u),\,j,k\in{\boldsymbol I}_n,$ and carries exactly $n$-multiple of the representation $E$.\,  
However, the stand point of \cite{PW27} is different from this point of view, as is explained below. 
\hfill
$\Box$
\vskip1.2em
  
For any unitary matrix representation $E(s)$, \lq Matrix\rq\, ${\boldsymbol A}(x)$ and \lq Zahl\rq\, $\chi(x)$ are defined as an analogy of Fourier transform\,:\; with $\chi(s):={\rm tr}(E(s))$, 
\begin{gather}
\label{2021-08-22-12}
{\boldsymbol A}(x):=\int x(s)E^*(s^{-1})\,ds =\int x(s)\overline{E}(s)\,ds\;\;(=\overline{E}(x)), 
\\
\nonumber
{\boldsymbol A}(x)=\big(\alpha_{ik}(x)\big),\quad \alpha_{ik}(x):=
\int x(s)\,\overline{e}_{ik}(s)\,ds\quad\mbox{\rm (Fourier-coefficients)}, 
\\
\alpha(x):={\rm tr}({\boldsymbol A}(x))=\int_{\mathfrak{G}} x(s)\overline{\chi}(s)\,ds. 
\end{gather}

In original notation in \cite{PW27}, $\widetilde{x}(s):=\overline{x(s^{-1})}$ and so, 
\;${\boldsymbol A}(\widetilde{x})={\boldsymbol A}(x)^*$. For $z=x\widetilde{x}$, ${\boldsymbol A}(z)={\boldsymbol A}(x){\boldsymbol A}(x)^*$ is positive Hermitian, and \;${\rm tr}({\boldsymbol A}(z))=\sum_{i,k}|\alpha_{ik}(x)|^2$. 
Fourier transform of $x$ is, in a sense, a decomposition of of the delta functional $S(x):= x({\small o})$ at the origin ${\small o}$, into a positive linear sum of Fourier coefficients of $x$. Applying it to Hermitian element $z=x\widetilde{x}$, we will get an $L^2$-type equality. 

To establish this equality (called here as \lq\,{\bf Die Vollst\"andichkeits der primitiven Darstellung}\,') is the main purpose of the paper. 

One step before the equality, we will get an inequality called 
\lq\,{\sc Bessel}sche Ungleichung\,\rq\/ given as follows: \;with $n=\dim E$,
and $E$ runs over representatives of all equivalence classes of irreducible unitary representations, 
\begin{gather}
\label{2021-08-23-1}
\sum_E n\,{\rm tr}\big({\boldsymbol A}(x\widetilde{x})\big)= 
\sum_E \dim E\sum_{i,k}|\alpha_{ik}(x)|^2 
\le \int|x(s)|^2ds. 
\end{gather} 

From this inequality, to arrive at the equality, the following two facts are important in my opinion.  
\vskip.7em

(1$^o$)\; The integral operator $I_K$ with continuous kernel $K$ is compact, 
that is, for any bounded sequence $f_1,f_2,\,\ldots$ in ${\cal H}$, we can find a 
subsequence $f_{i_1}, f_{i_2},\,\ldots$ such that the series $I_Kf_{i_k}$ 
converges as $k\to\infty$, or equivalently, the image of any bounded subset of 
${\cal H}$ is relatively compact.  

(2$^o$)\; For any $x\in C({\mathfrak{G}})$, the integral operator $I(x\widetilde{x})=I(x)I(x)^*$\, is Hermitian positive definite, and only has eigenvalues 
and eigenspaces. 
\vskip1em

Actually Peter and Weyl noted as \;
\begin{quotation}
{\it 
Die Vollst\"andigkeitsrelation gewinnen wir aus der Theorie der Eigenwerte und Eigenfunktionen von Kernen der besonderen Gestart $x(st^{-1})$. 
}
\end{quotation}

\section{Long-awaited general theory of invariant measures on groups}

Haar's result \cite{Haa33} is a long-awaited general theory of 
invariant measures on groups. Its main result is for locally compact 
groups, not necessarily Lie groups. 
\vskip1,2em 

{\bf Theorem 6.1.}\; 
{\it 
If a locally compact group is metrisable and separable, then there exists 
on it 
an outer measure which is right-invariant (or left-invariant).}
\vskip1.2em 

In its Introduction he referred to the above 
works of Hurwitz, Schur and Weyl, as we quote below from the top of 
Introduction:  
\begin{quotation}
{\bf 1.}\; 
Der Ausgangspunkt der {\sc Lie}schen Theorie der kontinuierlichen Gruppen, 
die sog.\;Infinitesimaltransformation, wird bekantlich mittels eines 
Differentiationsprozesses gewonnen; deshalb ist die {\sc Lie}sche Theorie 
in ihrer urspr\"unglichen Form nur auf solche Gruppen anwendbar, welche 
durch solche Gleichungen dargestellt sind, die die fraglichen 
Differentiationsbedingungen erf\"ullen.  
 Dieser Theorie steht eine andere gegen\"uber, die von Hurwitz in einer 
ber\"umten Arbeit angebahnt wurde, welche man treffend als eine 
Integrationstheorie 
der kontinuierlichen Gruppen bezeichnet hat; diese Theorie wurde 
insbesondere im letzten Jahrzehnt durch eine Reihe von wichtigen 
Arbeiten gef\"ord-
\\
ert, von denen wir hier nur die sch\"onen Arbeiten von 
Schur und Weyl erw\"ahnen.
\end{quotation}

It continues as 
\begin{quotation}
Es liegt daher der Gedanke nahe, die Frage zu unterzuchen, ob man 
$\cdots\cdots\cdots$\;\; 
Diese Frage ist offenbar damit gleichwertig, ob man \;{\it in der 
Gruppenmannigfaltigkeit einen Inhalts- bez.\;Ma\ss\/begriff einf\"uren kann, 
der invariant gegen\"uber den Transformationen der Gruppe ist}, \,d.\;h.\; 
der 
$\cdots\cdots$ (2 lines omitted) $\cdots\cdots$
\\
Unsere Untersuchungen gelten sogar f\"ur noch allgemeinere kontinuierliche 
Gruppen; wir werden im wesentlichen nur annehmen da\ss\/ die 
{\it Gruppenmannigfartigkeit metrisch, separabel und im Kleinen kompakt ist.}
\; $\cdots\cdots$
\end{quotation}

The paper \cite{Haa33} is separated into 15 Parts, numbered as 
{\bf 1}\,--\,{\bf 15}, and  
to introduce its contents, we make up Table of Contents as 
\vskip.7em

\qquad \S 1\; Der Inhalts. \hspace{33.8ex}Part 2--6\hspace{4.2ex} pp.148--155,  

\qquad \S 2\; Eigenshaften der Inhaltes.\hspace{20.5ex}Part 7--8\hspace{1.8ex}
\quad pp.155--160,

\qquad \S 3\; Das Analogon des {\sc Lebesgue}schen Ma\ss\/es. \hspace{1.8ex}
Part 9--12\hspace{3ex} 
pp.160--166, 

\qquad \S 4\; Anwendungen. \hspace{31ex}Part 13--15\hspace{1.9ex} pp.166--169.

\vskip.8em  
\noindent
In \S 4, there discussed three kinds of applications, and the third one is 
an extension of Theorem of Peter-Weyl on compact Lie groups to the case of 
a locally compact groups metrizable and separable. The top of Part 15 is 
\begin{quotation}
{\bf 15}. \;Ist die Gruppenmannigfaltigkeit ${\mathfrak{G}}$ kompakt, 
so kann mann 
ohne Schwierigkeiten die sh\"one Theorie von F.\,Peter und H.\,Weyl\footnote{
Mathematische Annnalen, Bd. 97, S.737--755.} 
\"uber die Darstellung der geschlossenen {\sc Lie}schen Gruppen auf den 
vorliegenden Fall \"ubertragen, da in diesen Untersuchungen lediglich nur 
der invariante Integrationsproze\ss\/ benutzt wird. Um dies kurz 
anzudeuten, $\cdots\cdots\cdots$\,. 
\end{quotation} 
As remarked before, the essential points of the proof are 
(1$^o$) and 
(2$^o$) at the end of \S 5.

\section{Two different proofs of uniqueness of Haar measure}

Neumann gave two different kind of proofs of the uniqueness of Haar measure 
in the papers \cite{Neu35} 
 and \cite{Neu36}.
\vskip.5em

{\bf \;7.1. \;The first paper \cite{Neu35}.}

 In its Part {\bf 2}, 
Neumann introduced, in place of 
Haar-Lebesgue type measure theory, another approach as follows. Let 
$C_{\boldsymbol R}(G)$ be the space of all real continuous functions on a 
compact group $G$. 
Consider real functional $M$ (called \lq\lq\,Mittel\,\rq\rq\,=\,mean) on $C_{\boldsymbol R}(G)$ which satisfies 
the conditions 1) to 7) below: 
\vskip.7em

\qquad 1)\;\; $M\big(\alpha f(x)\big)=\alpha M\big(f(x)\big)\quad(\alpha\in{\boldsymbol R})$.

\qquad 2)\;\; $M\big(f(x)+g(x)\big)=M\big(f(x)\big)+M\big(g(x)\big)$.

\qquad 3)\;\; If $f(x)\ge 0\;(x\in G)$, then \;$M\big(f(x)\big)\geq 0$.  

\qquad 4)\;\; If $f(x)=1\;(x\in G)$, then \;$M\big(f(x)\big)=1$.  

\qquad 5)\;\; If $M\big(f(xa)\big)=M\big(f(x)\big)\quad (a\in G)$.  

\qquad 6)\;\; If $M\big(f(ax)\big)=M\big(f(x)\big)\quad (a\in G)$.  

\qquad 7)\;\; If $M\big(f(x^{-1})\big)=M\big(f(x)\big)$. 
\\[1.2ex] 
Neumann proved that such an $M$ gives a Haar-Lebesgue type measure on $G$,  
saying as  
\begin{quotation}
\noindent
Mit der Hilfe 
eines solchen Mittels kann n\"amlich ein {\sc Haar-Lebesgue}sches Ma\ss\/ 
eingef\"uhrt werden, wie die folgenden \"Uberlegungen zeigen: 
$\cdots\cdots\cdots\cdots\cdots$
\end{quotation}

In Part {\bf 3}, he constructed such an $M$\, by managing \lq\,averages\,\rq\/ 
of right translations \;$f(xa)\;\,(=R(a)f(x))$,  
with very delicate technique, and get the so-called 
\lq\lq\,recht-Mittel\,\rq\rq. Then in Part {\bf 4}, he proved that the 
constructed 
$M$\, actually satisfies the conditions 1) to 7). Thus he proved 
the existence 
and the uniqueness of Haar measure, at the same time. 
The obtained $M$ is called \lq\lq\,Mittel stetiger Funktion\,\rq\rq.

\vskip1.2em 

{\bf Theorem 7.1.}\; {\it 
On a compact group $G$, without any set theoretical or topological 
 restriction, 
there exists a right-invariant Haar measure, 
unique up to a constant multiple. Moreover it is left-invariant and 
also invariant under the inversion $x\to x^{-1}$.
}\hfill
$\Box$
\vskip1.2em

{\bf \;7.2. \;The second paper \cite{Neu36}.}

The paper is divided into 11 Parts, and organized as follows: 
\vskip.8em

\qquad Introduction\hspace{27.4ex} Parts 1--2, 

\qquad Proof of the Theorem of uniqueness\qquad Parts 3--8, 

\qquad Consequences\hspace{26.3ex} Parts 9--10, 

\qquad Appendix\hspace{30.2ex} Part 11. 

\vskip.8em

In Part 2 of Introduction, Neumann remarked the following. 
 
The extension of the above 
process in case of compact groups 
to get \lq\lq\,Mittel stetiger Funktion\,\rq\rq\/ can be 
extended to non-compact case.  However the result is not an exterior 
measure at all, but the generalization of 
the integral mean for {\it almost periodic functions}\, of any group 
which is both left- and right-invariant. 

Then, as for the case of locally 
compact group, I quote from his text as  
\begin{quotation}
Thus, in the case of a general locally compact group an independent 
treatment of the problem of uniqueness is needed. 
This will be given in this paper: in fact, we shall prove the 
\vskip.6em 

{\bf Theorem of uniqueness.}\; 
{\it 
The left- as well as the right-invariant exterior measure in $G$ 
is unique.
\hfill
$\Box$
} 
\end{quotation}

\section{
Simple proof for extended 
existence-uniqueness theorem
}

{\bf 8.1.\; The first paper \cite{Kak36}.} 

In this paper, in connection to Haar's existence theorem 
of invariant measures on locally compact 
metrisable and separable groups, S.\;Kakutani treated 
the metrisability on topological groups. 
His results announced here is  
\begin{quotation}
{\it 
{\bf\it Satz.}\; Wenn die topologische Gruppe $G$ dem ersten 
Abz\"ahlbarkeitsaxiom gen\"ukt, dann kann man in $G$ eine 
metrik \,$\rho(x,y)$\, einf\"uhren, welche ausser den drei Distanzaxiomen 
noch der isomerischen Relation 
\begin{gather}
\label{2021-08-27-1}
\rho(zx,zy)=\rho(x,y)
\end{gather} 
gen\"ukt.
}
\end{quotation}
Thus, for a topological group, the metrisability is equivalent to the first 
countability axiom. 

To me, his proof here is very interesting and simple. 

\vskip1.2em

{\bf 8.2.\; The second paper \cite{Kak38}.} 

Let's quote from the 
top of Part {\bf 1} the following: 
\begin{quotation}
{\bf 1.}\; For a topological group $G$, which is locally compact and 
separable, the uniqueness of Haar's left-invariant measure is proved 
by J.\,v.\;Neumann. Although the method used by him is very 
interesting and powerful, his proof is somewhat long. The notion of 
right-zero-invariance is not necessary for the proof. In this paper 
we shall give a simplified proof. $\cdots\cdots\cdots$. 
\end{quotation}

Kakutani proposes moreover to treat the generalized case of $G\times S$, 
where $S$ is 
a topological space, and $G$\, acts transitively 
 on $S$\, as a group of homeomorphisms: 
\begin{quotation}
\noindent 
 $\cdots\cdots\cdots$\,.\;  
Since the separability plays no essential r\^ole in our proof, it can 
also be, by slight modifications applied to the case of a non-separable 
group (the case of a locally bicompact topological group, which is 
treated by A.\,Weil), and moreover we can prove, in the same manner, 
the theorem of the uniqueness of Haar's measure even for the case, 
when the field $G$ is no longer a topological group, that is, $G$ is simply 
a topological space $S$,  $\cdots\cdots\cdots$ 
\end{quotation}

A Borel set $E$ of $S$ is called {\it $\mu$-invariant}\; if 
$\mu(E\Delta \sigma E)=0\footnote{
$E_1\Delta E_2:=(E_1\setminus E_2)\bigsqcup\hspace{.3ex}(E_2\setminus E_1).$
} 
 \;(\forall \sigma \in G)$, and $G$ is called {\it ergodic on $S$}\; if 
for any invariant totally additive non-negative set function\footnote{
Not assumed that $\mu(U)>0$ for open set $U$.}
 $\mu$ on $S$ and for any $\mu$-invariant Borel set $E$ of $S$, either 
$\mu(E)=0$ or $\mu(S\setminus E)=0$.  In the special case where $S=G$, 
and $G$ acts 
on $S=G$\, through left (or right) translation, $\mu$-invariant set function on 
$S=G$\, is called {\it left} (or {\it right})-$\mu$-{\it invariant}\; and $G$ 
is called  {\it left} (or {\it right})-$\mu$-{\it ergodic}\; if $G$\, is 
ergodic on itself by respective translations.   

Let assume $G$\, be locally compact and {\it separable} (despite the previous 
statement).\footnote{Since the relation between uniqueness and 
ergodicity is more prominent in this case (by Kakutani).} 
In the case of $S=G$, he asserts as follows: 
\begin{quotation}
I.\; {\it If $G$ is left-ergodic, then the left-invariant measure of $G$ 
is unique (up to\\
\hspace*{5.6ex} a constant multiple).
}

II.\; $G$ is left-ergodic.
\end{quotation}

In Part {\bf 3}, the assertion I is proved, and 
here the conditions of separability and of local compactness are used. 
In Part {\bf 4}, the assertion II is proved. In Part {\bf 5}, the 
general case of non-separable locally compact groups $S=G$ is treated. 
Unfortunately, maybe because of the limitation of pages for 
Proc. Imp. Acad., detailed comments of the general 
case of $S$ is missing.

\section{
The definitive work on invariant measure theory on groups
}

The second chapter (\S\S 6--9) of this book treats invariant measures 
on topological groups and relative-invariant measures on quotient spaces. 
Table des Mati\`eres for Chapter 2 is 
\vskip.8em 

{\sc Chapitre}\; I\,I\qquad {\it La mesure de Haar}

\qquad \S\;6.\, Mesures et Int\'egrales,\; p.\;30.\quad\; \S\;7.\, Mesure 
de Haar,\; p.\;33.

\qquad\; \S\;8.\, Propri\'et\'es de la mesure de Haar,\; p.\;\;38.\qquad 
\S\;9.\,\; Mesures 

\quad dans les \'espaces homog\`enes,\; p.\;42. 
\vskip1em

In \S 6, Rodon measure is explained. 
For a locally compact topological space $X$, let ${\mathfrak{M}}_c(X)$ be the 
$\sigma$-algebra of subsets of $X$ 
generated by the set of all compact subsets, and a measure defined on 
 ${\mathfrak{M}}_c(X)$ is called {\bf Radon measure}. 
Radon measure can be defined 
by a real linear functional on the space $C_c(X)$ of real-valued 
continuous functions with compact supports on $X$. 
This way of treatment was used in \cite{Neu35} and later it is 
the way of Bourbaki to treat integrals.\footnote{
A.\,Weil himself was an important member of the founders of Bourbaki group.
} 
A real linear functional $\varphi$ on $C_c(X)$ is called {\it positive}\, 
if $\varphi(f)\ge 0$ for $f\ge 0,\;f\in C_c(X)$. The basic proposition here 
is 
\vskip1.2em
{\bf Proposition 9.1.}\; 
{\it For a positive real linear functional $\varphi$, there exists uniquely 
a Radon measure $\mu$ on ${\mathfrak{M}}_c(X)$ such that 
\begin{gather}
\label{2021-08-28-1}
\varphi(f)=\int_Xf(x)\,d\mu(x)\qquad(f\in C_c(X)).
\end{gather}
Conversely, for a Radon measure $\mu$ on $X$, consider the integral 
in the right hand side of (\ref{2021-08-28-1}), 
then it gives a positive real linear functional on $C_c(X)$. 
\hfill
$\Box$
}

\vskip1.2em

Let $X$ be a locally compact group $G$. Then a Radon measure $\mu$ is left-invariant if and only if so is the corresponding positive 
functional $\varphi$, that is, 
$\varphi(L_sf)=\varphi(f)$ for $s\in G$. Based on this fact, Weil proved 
\vskip1.2em

{\bf Theorem 9.2.}\; 
{\it 
Let $G$ be a locally compact group (not assumed any countability axiom). 
Then there exists a left-invariant positive functional $\varphi$, 
unique up to a constant multiple.
}
\vskip1em

The proof in \S 7 (pp.33--38) uses Zermelo's Axiom of Choice, 
according to Weil, in the form of Tychonoff's therem: \;
\lq\lq\,{\it A direct product of compact spaces is also compact}.\,''
\vskip1em

Furthermore, in Chapter 5 of the book, the so-called Peter-Weyl Theorem 
is proved in a complete form for general compact groups. 
\vskip1em

{\small {\bf Note.} 
The book \cite{Wei40}, known for its beautiful text, is 
finally translated into Japanese \cite{Sai15} by M. Saito, 
waited from long ago. 
}

\section{
Equivalency between Borel measure and Baire measure
}

The paper \cite{KK44} consists of short Introduction and four sections. 
Let's quote the top part of Introduction to see 
the problem studied here:
\begin{quotation}
Zur Definition des {\sc Haar}schen Masses $m$ in einer lokal bikompakten, 
nicht separabeln Gruppe gibt es zwei M\"oglichkeiten. Nach der 
ersten gew\"onlichen Definition wird $m$ zun\"achst f\"ur 
alle {\sc Borel}schen Mengen erkl\"art und dann zum vollst\"andigen Mass 
vervollst\"andigt\,;\; 
nach der zweiten wird dagegen $m$ zun\"achst nur f\"ur die M\"angen 
mit {\sc Baire}schen charakteristische Funktionen\footnote{
Let $G$ be a topological group. By definition, the space ${\cal B}(G)$ of 
all Baire functions on $G$ is the completion of the space $C(G)$ of 
real-valued continuous functions with respect to the repetitions of 
oparation 
\lq\lq\,for a series of functions $f_1,f_2,f_3,\cdots$ take their 
pointwise limit $f=\lim_{n\to\infty}f_n$\,\rq\rq
}
 --- wir wollen 
solche M\"ange {\sc Baire}sche nennen --- definiert und dann 
vervollst\"andigt. Sind nun diese zwei Definitionen \"aquivalen\;?\; 
In der vorliegenden Note soll diese Frage bejahend beantworted werden. 
$\cdots\cdots\cdots\cdots$
\end{quotation}

\vskip.3em
The principal result of this paper is the following one for 
left-invariant measures on locally compact groups\,: 
\vskip1.2em 

{\bf Theorem 10.1.}\; 
{\it 
The first method and the second method are mutually equivalent, that is, 
they give the same completed measure. 
}

\vskip2em

{\small  
Takeshi HIRAI, \;\;22-8 Nakazaichi-Cho, Iwakura, Sakyo-Ku, Kyoto 606-0027, Japan; 

hirai.takeshi.24e@st.kyoto-u.ac.jp
} 

\end{document}